\documentclass[review,hidelinks,onefignum,onetabnum]{siamart250211}
\usepackage{amssymb}
\renewcommand{\footercopyright}{}

\ifpdf
\hypersetup{
  pdftitle={Corruption via Mean Field Games},
  pdfauthor={Michael V. Klibanov, Mikhail
Yu. Kokurin, and Kirill V. Golubnichiy}
}
\fi
\begin{document}
\nolinenumbers

\title{Corruption via Mean Field Games}
\author{Michael V. Klibanov \thanks{
Department of Mathematics and Statistics, University of North Carolina at
Charlotte, Charlotte, NC, 28223, USA, mklibanv@charlotte.edu}, \and Mikhail
Yu. Kokurin \thanks{%
Department of Physics and Mathematics, Mary State Universiry, Yoshkar--Ola,
424000, Russian Federation, kokurinm@yandex.ru} \ and \and Kirill V.
Golubnichiy\thanks{%
Department of Mathematics and Statistics, Texas Tech University, Lubbock, TX
79409, USA, kgolubni@ttu.edu},}
\date{}
\maketitle

\begin{abstract}
A new mathematical model governing the development of a corrupted hierarchy
is derived. This model is based on the Mean Field Games theory. A
retrospective problem for that model is considered. From the applied
standpoint, this problem amounts to figuring out the past activity of the
corrupted hierarchy using the present data for this community. Three new
Carleman estimates are derived. These estimates lead to H\"{o}lder stability
estimates and uniqueness results for both that retrospective problem and its
generalized version. H\"{o}lder stability estimates characterize the
dependence of the error in the solution of the retrospective problem from
the error in the input data.
\end{abstract}

\textbf{Key words}. mathematical model, corrupted hierarchy, mean field
games, retrospective problem, three new Carleman estimates, H\"{o}lder
stability estimates, uniqueness of the solution.

\textbf{\ MSC codes}. 91A16, 35R30

\section{Introduction}

\label{sec:1}

In this paper, we introduce a new mathematical model of the development of a
corrupted hierarchy. This model is formulated in terms of the Mean Field
Games (MFG) theory. Next, we derive a version of the Mean Field Games System
(MFGS) of two coupled nonlinear parabolic PDEs with the opposite directions
of time. From the applied standpoint, an interesting question is: \emph{%
Given the present stage of a corrupted hierarchy, what was its historical
development?} We formulate this question as the retrospective problem for
the MFGS. Another argument in favor of an interest of the retrospective
problem for our mathematical model is that statistical \ data are
insufficient sometimes for figuring out the initial distribution $m\left(
x,0\right) $ of the density function $m\left( x,t\right) $. On the other
hand, the function $m\left( x,0\right) $ is conventionally used in the MFGS 
\cite{A}. To be more specific, we consider the case when the terminal
conditions 
\begin{equation}
u\left( x,T\right) ,\text{ }m\left( x,T\right) \text{ }  \label{1.1}
\end{equation}%
are known for both the minimal average cost function $u\left( x,t\right) $
and the density function $m\left( x\mathbf{,}t\right) ,$ where $T>0$ is the
final/present moment of time. The minimal average cost function $u\left(
x,t\right) $ is an analog of the value function in the MFG theory \cite{A}.

Here and below $x=\left( x_{1},...,x_{n}\right) \in \mathbb{R}^{n}$ denotes
the vector of spatial variables and $t\in \left( 0,T\right) $ is time. In
the particular case of our model $n=2,$ and we denote below $\mathbf{x}%
=\left( x,y\right) .$ Although another retrospective problem (not linked to
our mathematical model) for the MFGS was considered in \cite{MFG1,MFGbook},
it was assumed in these references that the knowledge of functions (\ref{1.1}%
) is complemented by the knowledge of the function 
\begin{equation}
m\left( x,0\right) .  \label{1.2}
\end{equation}%
The additional condition (\ref{1.2}) has resulted in the Lipschitz stability
estimate in \cite{MFG1,MFGbook}.

Any input data are given with an error, so as the input data for our
retrospective problem. Hence, it is important to estimate how that error
influences the error in the solution of our retrospective problem. We
address this question via derivations of H\"{o}lder stability estimates for
both our particular retrospective problem and its generalized version. In
particular, these estimates imply uniqueness of the solution of each of
these two problems. Since the above mentioned Lipschitz stability estimate
of \cite{MFG1,MFGbook} has led to a globally convergent numerical method for
that case \cite{MFG7,MFGbook}, then we believe that the H\"{o}lder stability
estimates of this paper might eventually lead to a globally convergent
numerical method for the retrospective problem we consider here. This, in
turn would allow one to conduct numerical studies of our model.

The MFG theory studies the behavior of infinitely many rationally acting
agents. This theory was introduced in 2006 in the seminal works of Huang,
Caines, and Malham\'{e} \cite{Huang2,Huang1} and Lasry and Lions \cite%
{LL1,LL2}. The MFG theory is broadly applicable to descriptions of many
complex social phenomena. Among those applications, we mention, e.g. finance 
\cite{A,Trusov}, sociology \cite{Bauso}, election dynamics \cite{Chow}. We
refer to the book \cite[chapter 6]{Kol} for more applications. In particular
this book considers corruption modeling via the MFG theory. However, our
MFG-based model of the corrupted hierarchy is significantly different in
many aspects from the one of \cite{Kol}. For example, the model of \cite{Kol}
is a stationary one. The retrospective problem cannot be considered for a
stationary model, which is unlike our case of the time dependent model. In
addition, our model is continuous, whereas the model of \cite{Kol} is
discrete.

The key to our H\"{o}lder stability estimates are three new Carleman
estimates for the MFGS, which we derive here. Carleman estimates were first
introduced in the MFG theory in \cite{MFG1}. As mentioned above, this
technique allowed to prove Lipschitz stability estimate for the above
outlined version of the retrospective problem when all three functions in (%
\ref{1.1}), (\ref{1.2}) are known. Later the tool of Carleman estimates was
applied to obtain both H\"{o}lder and Lipschitz stability estimates for
various problems for the MFGS \cite{MFG2,MFG6,MFG3,MFG4,Liao}. We also refer
to the recently published book \cite{MFGbook} on this subject. In addition,
this tool allows one to construct globally convergent numerical methods for
various problems for the MFGS, including coefficient inverse problems \cite%
{MFG7,MFGIPI,MFGbook}.

H\"{o}lder stability estimate for an analog of the retrospective problem of
this paper was obtained in \cite{MFG4}. However, the principal parts of the
PDE operators in \cite{MFG4} are $\partial _{t}\pm a\Delta ,$ where $a>0$ is
a number, i.e. this is the case of constant coefficients in the principal
part of the parabolic operators. Unlike this, our mathematical model
requires that the principal parts of parabolic operators of the MFGS should
be $\partial _{t}\pm L,$ where $L$ is an elliptic operator of the second
order with variable coefficients. Therefore, it is necessary to prove here
new Carleman estimates for the operators $\partial _{t}\pm L$. Another new
element of this paper is that while the zero Neumann boundary condition at
the whole boundary is used in \cite{MFG4} for both functions $u\left(
x,t\right) $ and $m\left( x,t\right) ,$ in our case each of these functions
has the Dirichlet boundary condition on a part of the boundary and the
Neumann boundary condition on the rest of the boundary.

In all above cited publications about applications of Carleman estimates to
the MFGS, the case of a single measurement input data is considered. We
refer to a series of recent publications \cite{Ding,Liu1,Liu2,Ren1,Ren2},
where inverse problems with multiple measurements for the MFGS are studied.

All functions considered below are real valued ones. In section 2 we provide
an informal description of our mathematical model. In section 3 we discuss
possible inverse problems for our model. The MFGS for our case is derived in
section 4, and then our retrospective problem is formulated in that section.
In section 5 we prove two new Carleman estimates mentioned above. The H\"{o}%
lder stability estimate for \ generalized retrospective problem in is
obtained in section 6. In section 7, the estimate of section 6 is specified
for the retrospective problem of section 4.

\section{Informal Description of the Model}

\label{sec:2}

Consider a community consisting of an infinite number of homogeneous agents,
where the current state of each agent at any time $t\in \lbrack 0,T]$ is
characterized by two parameters: the relative degree of corruption $x$ and
the relative position $y$ in the organizational hierarchy. We assume that $%
x\in \lbrack 0,1]$ and $y\in \lbrack 0,1]$, where $x=0$ represents a
complete lack of corruption and $y=0$ means the lowest level in the
hierarchy.

At any moment in time, an agent has control actions $\alpha ,\beta
,u=(\alpha ,\beta )\in U$, where $\alpha \geq 0$ represents efforts to
advance in the illegal (corruption) hierarchy, and $\beta \geq 0$ denotes
efforts to ascend the organizational ladder. Here, $U\subset \mathbb{R}^{2}$
represents the set of admissible controls. Negative values of $\alpha $ and $%
\beta $ correspond to actions aimed at moving towards honest behavior (e.g.,
partial or complete rejection of corruption) or voluntarily stepping down to
a lower position in the organizational hierarchy. Choosing $\alpha =\beta =0$
over a certain period of time indicates that the agent does not undertake
any active measures to change the status that agent has.

The financial income of an agent at any given time $t$ at the state $(x,y)$
is denoted as $c(x,y,t)$. This income can be written as: 
\[
c(x,y,t)=p(y,t)+q(x,y,t), 
\]%
where $p(y,t)$ represents lawful salary, and $q(x,y,t)$ denotes the unlawful
income from the corruption activities of the corruption degree $x$. The
functions $p$ and $q$ are naturally assumed to be increasing with respect to
their arguments. The dependence on $t$ is due to such factors as salary
indexing and inflation.

The vector of control actions $(\alpha ,\beta )$ generates the cost $%
h(\alpha ,\beta )$ per unit time. In the simplest case, the cost can be
modeled as: 
\begin{equation}
h(\alpha ,\beta )=\frac{1}{2}a_{0}(\mathbf{x})\alpha ^{2}+\frac{1}{2}b_{0}(%
\mathbf{x})\beta ^{2},\quad a_{0}(\mathbf{x}),b_{0}(\mathbf{x})>0.
\label{2.01}
\end{equation}

More complex functions can also be considered to reflect the asymmetry of
the costs associated with increasing or decreasing $x$ and $y$. For example,
increasing $x$ might involve financial contributions to the corrupted
networks. When $x$ decreases, the agent pays a compensation to the corrupted
community, such as penalties for either interrupting or narrowing schemes of
their enrichment or for both of these. When $y$ increases, the agent pays a
financial cost to validate the higher status of this individual (e.g.,
purchasing more expensive goods). Conversely, when $y$ decreases, the cost
may involve organizing a transfer to a desired lower position. Both types of
transfers may depend on $x$, which is reflected in the potential dependence
of the values $a_{0}$ and $b_{0}$ on $x.$

The rationale behind efforts of transitions to lower positions is to
minimize the attention of supervisory authorities, which tend to scrutinize
more closely individuals at the higher levels of the hierarchy. The function 
$h(\alpha ,\beta )$ can also describe the intellectual and emotional efforts
made by an agent to implement the controls $(\alpha ,\beta ).$ In any case,
this function should attain its minimum value at the point $(0,0)$ as a
function of $(\alpha ,\beta ).$

When making decisions regarding the choice of controls $(\alpha ,\beta )$,
an agent considers not only the desire to minimize the total costs over the
operational time interval $[0,T]$, but also the current state of the entire
community. At each moment of time $t$, this state is described by the
density $m(\mathbf{x},t)$, representing the distribution of agents across
states 
\begin{equation}
\left. \mathbf{x}=(x,y)\in \Omega =[0,1]\times \lbrack 0,1].\text{ }\right.
\label{2.2}
\end{equation}

The quantity $m(x, y, t) \geq 0 $ is proportional to the number of agents
that are in the state $(x, y)$ at time $t$. We adopt the standard assumption
of the MFG theory that the community bases constructs its controls depending
on this density. A typical objective for each individual agent is to
minimize the total cost of this person while behaving as all other members
of this community.

To formalize this behavior, we introduce the functions $\overline{m}%
^{(y)}(x,t)$ and $\overline{m}^{(x)}(y,t)$

\begin{equation}
\left. 
\begin{array}{c}
\overline{m}^{(y)}(x,t)=\left( \varepsilon
+\int\limits_{0}^{1}m(x,y,t)\,dy\right)
^{-1}\int\limits_{0}^{1}y\,m(x,y,t)\,dy,\quad \\ 
\overline{m}^{(x)}(y,t)=\left( \varepsilon
+\int\limits_{0}^{1}m(x,y,t)\,dx\right)
^{-1}\int\limits_{0}^{1}x\,m(x,y,t)\,dx,\quad \varepsilon >0%
\end{array}%
\right.  \label{2.3}
\end{equation}

\noindent In the case $\varepsilon =0$, formulas (\ref{2.3}) describe the
average densities in terms of the levels of corruption and positions within
the hierarchy. To avoid technical difficulties associated with the
degenerate case when $m(x,y,t)=0$ on an interval of $x$ or $y$ lying inside $%
\Omega $, it is convenient for us to assume that $\varepsilon >0.$ Thus, $%
\overline{m}^{(y)}(x,t)$ represents the average position within the
hierarchy held by agents with a corruption level $x$, and $\overline{m}%
^{(x)}(y,t)$ represents the average corruption level of agents at a
hierarchical position $y$.

In addition to financial indicators, an agent may also want to minimize the
deviation of the current state $(x,y)$ of this individual from the
corresponding average $(\overline{m}^{(y)}(x,t),\overline{m}^{(x)}(y,t))$ at
any given moment in time.

In the case under consideration, $\varepsilon > 0$ implies that an arbitrary
agent is oriented toward slightly underestimated average state values,
averaged over the ensemble of agents with fixed characteristics $y$ or $x$,
respectively. The choice $\varepsilon = 0$ corresponds to targeting the
exact averaged states.

The above deviation is measured by the following function 
\begin{equation}
g(x-\overline{m}^{(y)}(x,t),y-\overline{m}^{(x)}(y,t)),  \label{2.4}
\end{equation}%
where the function $g\left( \mathbf{x}\right) $ is smooth, 
\begin{equation}
g\in C^{1}\left( \mathbb{R}^{2}\right) .  \label{2.5}
\end{equation}%
This function should achieve its minimal value at the point $(0,0)$. In the
simplest case, the function $g$ can be expressed as: 
\[
g(x,y)=\frac{1}{2}a_{1}x^{2}+\frac{1}{2}b_{1}y^{2},\quad a_{1},b_{1}>0.
\]%
In the general case, the motion of the agent in the phase space $\Omega $
may terminate before the previously fixed time $T$. This occurs when the
point $(x,y)$, which describes the state of that agent, reaches the absorbing
part of the boundary of $\Omega $. Thus, the agent's dynamics takes place
for $t\in \lbrack 0,\hat{T}]$, where $\hat{T}\leq T$. The functional
describing the total financial and intellectual costs of an agent over the
time interval $[0,\hat{T}]$ also includes the term $\Psi (x_{\hat{T}},y_{%
\hat{T}})$, which depends on the state $(x_{\hat{T}},y_{\hat{T}})$ at the
final moment of time. If $\Psi (x_{\hat{T}},y_{\hat{T}})>0,$ then this term
denotes the profit of an agent in the case when the random walk of this
agent ends at the final moment of time $t=\hat{T}.$ In addition to the
standard \textquotedblleft severance pay", which is proportional to the
salary $p(y_{\hat{T}},\hat{T})$ at the final moment of time $t=\hat{T}$,
this term may also include, for example, a penalty paid in the case of an
exposure of this agent to authorities. In the case of a penalty $\Psi (x_{%
\hat{T}},y_{\hat{T}})<0.$ In what follows, we assume that the random walk is
terminated when an agent reaches the absorbing part $\left\{ x=1\right\} $
of the boundary; see below. Thus, we assume that the most corrupt agents are
immediately removed from the group of agents, for example, as a result of an
exposure to authorities. 

The movement of an agent in the phase space $\Omega $ is influenced by both
deliberate controls and random effects. The controlled system describing the
agent's dynamics is the following system of two stochastic differential
equations: 
\begin{equation}
\left. 
\begin{array}{c}
dx_{t}=\alpha \varphi _{1}\left( x_{t},y_{t}\right) dt+\sigma _{1}\left(
x_{t},y_{t}\right) dW_{1t}, \\ 
dy_{t}=\beta \varphi _{2}\left( x_{t},y_{t}\right) dt+\sigma _{2}\left(
x_{t},y_{t}\right) dW_{2t}.%
\end{array}%
\right.  \label{2.6}
\end{equation}
Here, $(x_{t},y_{t})$ represents the agent's position in the phase space at
time $t$, and $W_{1t}$, $W_{1t}$ are two independent standard Wiener
processes (one-dimensional Brownian motions). The terms $\sigma _{1}(x,y)$, $%
\sigma _{2}(x,y)>0$ denote the volatilities of these processes. It is
reasonable to assume that the functions $\sigma _{1}$ and $\sigma _{2}$
decrease with respect to each of their two arguments. If $\alpha =\beta =0$,
then the dynamics of an agent essentially reduces to two-dimensional
Brownian motion in $\Omega $, meaning that movements of agents in the $x$
and $y$ directions are purely random.

The factors $\varphi _{1}(\mathbf{x})$ and $\varphi _{2}(\mathbf{x})$ play a
key role in the model, as they describe the amplification (or attenuation)
of the control effects $\alpha $ and $\beta $ on the agent's speed in the $x$
and $y$ directions. Since the phase space $\Omega $ is bounded, then the
following formulas for functions $\varphi _{1},\varphi _{2}$ are considered
reasonable ones 
\begin{equation}
\varphi _{1}(x,y)=ax\left[ (1-x)+p_{1}y\right] ,\quad \varphi _{2}(x,y)=by%
\left[ (1-y)+p_{2}x\right]  \label{2.7}
\end{equation}%
and also%
\begin{equation}
\varphi _{1}(x,y)=ax\left[ (1-x)-p_{1}y\right] ,\quad \varphi _{2}(x,y)=by%
\left[ (1-y)-p_{2}x\right] ,  \label{2.8}
\end{equation}%
with some numbers $a,b,p_{1},p_{2}>0$.

For example, if $\sigma _{2}\equiv 0$, then it follows from equations (\ref%
{2.6})-(\ref{2.8}) that, in the absence of random disturbances and
corruption ($x=0$), a deterministic monotonic career growth takes place.
This growth is governed by the equation: 
\begin{equation}
\dot{y}_{t}=\beta by_{t}(1-y_{t}).  \label{2.9}
\end{equation}%
Equation (\ref{2.9}) has two stationary solutions: $y_{t}\equiv 0$ and $%
y_{t}\equiv 1$. If $\beta >0,$ then for any initial value $%
y_{t_{0}}=y^{(0)}\in (0,1)$ the solution $y_{t}$ of this equation is a
monotonically increasing function, and 
\begin{equation}
\lim_{t\rightarrow \infty }y_{t}=1,  \label{2.10}
\end{equation}%
i.e., under normal conditions, reaching the upper levels of the hierarchy is
possible only as $t\rightarrow \infty $, which corresponds to the reality.

Let $x>0$. Then equation (\ref{2.7}) implies that the corruption component
accelerates its upward movement along the $y-$axis as well as the downward
movement along this axis. In contrast, in the model described by equation (%
\ref{2.8}), corruption slows down the career shifts in both directions.

In the case of equation (\ref{2.7}), a high position $y$ in the hierarchy
accelerates the corruption process. In the model described by equation (\ref%
{2.8}), the opposite is true. Of course, other combinations of signs for the
second terms in the square brackets in equations (\ref{2.7}) and (\ref{2.8})
are possible. It would be of an interest to analyze from this perspective
typical examples of corrupted bureaucratic structures.

The controlled dynamics of the population of agents is described by the
functions $u(\mathbf{x},t)$ and $m(\mathbf{x},t)$, where $\mathbf{x}\in
\Omega $ and $t\in \lbrack 0,T]$. Here, $u(\mathbf{x},t)$ is the minimal
average cost for an agent starting at the position $\mathbf{x}$ at the
moment of time $t$ over the operational interval $[t,T]$. The function $m(%
\mathbf{x},t)$ denotes the distribution of agents across states $\mathbf{x}%
=(x,y)$ at the moment of time $t$.

These functions satisfy a nonlinear system of coupled integral differential
parabolic equations, which is a specific version of MFGS. In the
conventional formulation of these equations, the initial condition is given
for the function $m$, $m(\mathbf{x},0)=m_{0}(\mathbf{x})$, and the terminal
condition is given for the function $u$, as $u(\mathbf{x},T)=u_{T}(\mathbf{x}%
)$, where $\mathbf{x}\in \Omega $. The function $m_{0}(\mathbf{x})$
represents the degree of corruption across various levels of the hierarchy
at the initial moment of time, while $u_{T}(\mathbf{x})$ is denoted as: 
\[
u_{T}(\mathbf{x})=\Psi (\mathbf{x}). 
\]

A crucial aspect is the assignment of boundary conditions on $\Gamma
=\partial \Omega $. It is assumed that the randomly controlled trajectory $%
(x_{t},y_{t})$, described by equation (\ref{2.6}), is absorbed at the
absorbing portion $\Gamma _{0}$ of the boundary of the square $\Omega $ in (%
\ref{2.2}), 
\begin{equation}
\Gamma _{0}=\{(x,y)\in \partial \Omega :x=1\},  \label{2.11}
\end{equation}%
while at the remaining part 
\begin{equation}
\Gamma _{1}=\partial \Omega \setminus \Gamma _{0}  \label{2.12}
\end{equation}%
it is reflected off the boundary. Absorption signifies the removal of a
corrupt agent from the community. from the system. Reflection represents the
presence of managerial and societal mechanisms which prevent shifts below $%
y=0$ and above $y=1$ (since such positions do not exist), as well as keeping
agents to the right of the line $\left\{ x=0\right\} $.

Optimal controls $\alpha _{t}=\alpha _{t}(\mathbf{x})$ and $\beta _{t}=\beta
_{t}(\mathbf{x})$ are constructed based on feedback schemes and are
determined as solutions to of an initial boundary value problem for the MFGS
involving functions $u\left( \mathbf{x},t\right) $ and $m\left( \mathbf{x}%
,t\right) $.

\section{Possible Inverse Problems and Their Purpose}

\label{sec:3}

In the conventional formulation of an initial boundary value problem one
assumes that the terminal condition $u(\mathbf{x},T)$ for the function $u$
is known and the initial condition $m(\mathbf{x},0)$ for the function $m$ is
also known \cite{A}. In addition, a boundary condition is known for each of
these functions. However, uniqueness theorems for this case are proven only
under quite restrictive conditions \cite{A,LL2}.

We consider here a retrospective inverse problem. In this case, the
functions $u(\mathbf{x},T)$ and $m(\mathbf{x},T)$ are given (see (\ref{2.2}%
)), and the goal is to recover the initial distributions $u(\mathbf{x},0)$
and $m(\mathbf{x},0).$ Unlike the standard formulation, this approach can be
motivated by the lack of a detailed statistics for $m(\mathbf{x},0)$.

In this context, the interest may not lie solely in $m(\mathbf{x},0)$ but
also in the control functions $\alpha _{t}=\alpha _{t}(\mathbf{x})$ and $%
\beta _{t}=\beta _{t}(\mathbf{x})$, which characterize the psychological
part of the collective consciousness about the corruption.

Another formulation of an inverse problem involves the reconstruction of
certain $\mathbf{x}-$dependent coefficients of the mathematical model given
below. It is likely that the government would be particularly interested in
the functions $\varphi _{1}\left( \mathbf{x}\right) $ and $\varphi
_{2}\left( \mathbf{x}\right) $ in equations (\ref{2.7}), (\ref{2.8}).
However, we focus in this paper only on the retrospective problem.

\section{The MFG System and the Statement of the Retrospective Problem}

\label{sec:4}

Let $\tau $ denotes the time when the trajectory $(x_{t},y_{t})$, described
by the stochastic differential equations (\ref{2.6}), reaches the absorbing
part $\Gamma _{0}$ of the boundary $\Gamma $, and let $\widehat{T}=\min
\{\tau ,T\}$. The population of agents solves the problem of the
minimization of the mathematical expectation $\mathbb{E}$ of the total cost
over the time interval $[0,\widehat{T}]$:%
\begin{equation}
\left. 
\begin{array}{c}
\min_{\left( \alpha _{t},\beta _{t}\right) \in U,t\in \left[ 0,\widehat{T}%
\right] } \\ 
\mathbb{E}\left\{ 
\begin{array}{c}
\int\limits_{0}^{\widehat{T}}\left[ -c\left( x_{t},y_{t},t\right) +h\left(
\alpha _{t},\beta _{t}\right) +g(x_{t}-\overline{m}^{(y)}(x_{t},t),y_{t}-%
\overline{m}^{(x)}(y_{t},t)\right] dt+ \\ 
+\Psi \left( x_{\widehat{T}},y_{\widehat{T}}\right)%
\end{array}%
\right\} .%
\end{array}%
\right.  \label{4.1}
\end{equation}

Equations (\ref{2.6}) imply that, the Hamiltonian of the controlled system
is:%
\begin{equation}
\left. 
\begin{array}{c}
H(x,y,t,m,p,\alpha ,\beta )= \\ 
=(\alpha \varphi _{1}(x,y)p_{1}+\beta \varphi
_{2}(x,y)p_{2}-c(x,y,t)+h(\alpha ,\beta )+ \\ 
+g\left( x-\overline{m}^{(y)}(x,t),y-\overline{m}^{(x)}(y,t)\right) ,\text{ }%
p=\left( p_{1},p_{2}\right) .%
\end{array}%
\right.  \label{4.2}
\end{equation}
The controls $(\alpha _{t},\beta _{t})=(\alpha _{t}(x_{t},y_{t}),\beta
_{t}(x_{t},y_{t}))$ are determined using a feedback scheme via solution of
the following minimization problem: 
\[
\min_{(\alpha ,\beta )\in U}H(x,y,t,m,p,\alpha ,\beta ). 
\]%
Let a solution of this problem be: 
\[
(\alpha ^{\ast },\beta ^{\ast })=(\alpha ^{\ast }(x,y,t,m,p),\beta ^{\ast
}(t,x,y,t,m,p)). 
\]%
Then the resulting optimal controls are:%
\begin{equation}
\left. 
\begin{array}{c}
\alpha _{t}=\alpha ^{\ast }(x_{t},y_{t},t,m(x_{t},y_{t},t),\nabla
u(x_{t},y_{t},t)), \\ 
\beta _{t}=\beta ^{\ast }(x_{t},y_{t},t,m(t,x_{t},y_{t}),\nabla
u(x_{t},y_{t},t)).%
\end{array}%
\right.  \label{4.3}
\end{equation}

Assume that the function $h$ has the form (\ref{2.01}) and that $U=\mathbb{R}%
^{2}$. Then the optimal controls are simplified: 
\[
\alpha ^{\ast }=-\frac{\varphi _{1}(\mathbf{x})}{a_{0}(\mathbf{x})}%
p_{1},\quad \beta ^{\ast }=-\frac{\varphi _{2}(\mathbf{x})}{b_{0}(\mathbf{x})%
}p_{2}. 
\]

Thus, the optimal controls are: 
\[
\alpha _{t}=-\frac{\varphi _{1}(x_{t},y_{t})}{a_{0}(x_{t},y_{t})}%
u_{x}(x_{t},y_{t},t),\quad \beta _{t}=-\frac{\varphi _{2}(x_{t},y_{t})}{%
b_{0}(x_{t},y_{t})}u_{y}(x_{t},y_{t},t). 
\]

Denote 
\begin{equation}
\left. 
\begin{array}{c}
Q_{T}=\Omega \times \left( 0,T\right) ,S_{T}=\partial \Omega , \\ 
\Gamma _{0,T}=\Gamma _{0}\times \left( 0,T\right) ,\Gamma _{1,T}=\Gamma
_{1}\times \left( 0,T\right) .%
\end{array}%
\right.  \label{4.04}
\end{equation}

Following the well-known scheme (see, e.g., \cite[pages 139, 327-328]{Carm}%
), we obtain the MFGS system with respect to two unknown functions $u=u(%
\mathbf{x},t)$ and $m=m(\mathbf{x},t).$

The first equation of the MFGS for functions $u(\mathbf{x},t)$ and $m(%
\mathbf{x},t)$ is:%
\[
u_{t}+\frac{\sigma _{1}^{2}(\mathbf{x})}{2}u_{xx}+\frac{\sigma _{2}^{2}(%
\mathbf{x})}{2}u_{yy}- 
\]%
\begin{equation}
-\left( \frac{\varphi _{1}^{2}(\mathbf{x})}{2a_{0}(\mathbf{x})}u_{x}^{2}+%
\frac{\varphi _{2}^{2}(\mathbf{x})}{2b_{0}(\mathbf{x})}u_{y}^{2}\right) -
\label{4.4}
\end{equation}%
\[
-c(\mathbf{x},t)+g(x-\overline{m}^{(y)}(x,t),y-\overline{m}^{(x)}(y,t))=0,%
\text{ }\left( \mathbf{x},t\right) \in Q_{T}. 
\]%
The second equation is:%
\[
m_{t}-\frac{1}{2}\left( \sigma _{1}^{2}(\mathbf{x})m\right) _{xx}-\frac{1}{2}%
\left( \sigma _{2}^{2}(\mathbf{x})m\right) _{yy}+ 
\]%
\begin{equation}
+\left( \frac{\varphi _{1}^{2}(\mathbf{x})}{a_{0}(\mathbf{x})}mu_{x}\right)
_{x}+\left( \frac{\varphi _{2}^{2}(\mathbf{x})}{b_{0}(\mathbf{x})}%
mu_{y}\right) _{y}=0,\text{ }\left( \mathbf{x},t\right) \in Q_{T}.
\label{4.5}
\end{equation}

We note that by (\ref{2.3}), the term $g(x-\overline{m}^{(y)}(x,t),y-%
\overline{m}^{(x)}(y,t))$ in (\ref{4.4}) has the form:%
\begin{equation}
\left. 
\begin{array}{c}
g(x-\overline{m}^{(y)}(x,t),y-\overline{m}^{(x)}(y,t))= \\ 
=g\left( 
\begin{array}{c}
x-\left( \varepsilon +\int\limits_{0}^{1}m(x,y,t)\,dy\right)
^{-1}\int\limits_{0}^{1}y\,m(x,y,t)\,dy, \\ 
y-\left( \varepsilon +\int\limits_{0}^{1}m(x,y,t)\,dx\right)
^{-1}\int\limits_{0}^{1}x\,m(x,y,t)\,dx,%
\end{array}%
\right) ,\varepsilon >0.%
\end{array}%
\right.  \label{4.50}
\end{equation}%
We also recall that 
\begin{equation}
m(x,y,t)\geq 0.  \label{4.51}
\end{equation}%
Hence, by (\ref{2.4}), (\ref{2.5}), (\ref{4.50}) and (\ref{4.51})%
\begin{equation}
\left. 
\begin{array}{c}
\left\vert g(x-\overline{m}_{1}^{(y)}(x,t),y-\overline{m}%
_{1}^{(x)}(y,t))-g(x-\overline{m}_{2}^{(y)}(x,t),y-\overline{m}%
_{2}^{(x)}(y,t))\right\vert \leq \\ 
\leq B\left( \left\vert \left( \overline{m}_{1}^{(y)}-\overline{m}%
_{2}^{(y)}\right) (x,t)\right\vert +\left\vert \left( \overline{m}_{1}^{(x)}-%
\overline{m}_{2}^{(x)}\right) (y,t)\right\vert \right) ,\text{ }x,y\in
\Omega ,t\in \left( 0,T\right) , \\ 
\forall m_{1},m_{2}\in C\left( \overline{Q}_{T}\right) , \\ 
B=B\left( \max \left( \left\Vert m_{1}\right\Vert _{C\left( \overline{Q}%
_{T}\right) },\left\Vert m_{2}\right\Vert _{C\left( \overline{Q}_{T}\right)
}\right) ,\max_{\mathbb{R}^{2}}\left( \left\vert g\right\vert ,\left\vert
\nabla g\right\vert \right) ,\varepsilon \right) >0,%
\end{array}%
\right.  \label{4.52}
\end{equation}%
where the number $B$ depends only on listed parameters.

Recall that parts $\Gamma _{0}$ and $\Gamma _{1}$ of the boundary $\partial
\Omega $ of the domain $\Omega $ are defined in (\ref{2.11}) and (\ref{2.12}%
), where $\Gamma _{0}$ is the absorbing part of the boundary. Let $%
\partial_\nu$ be the outward normal derivative on $\Gamma_1$. The
boundary conditions for the function $u$ are: 
\begin{equation}
u|_{\Gamma _{0},_{T}}=u\mid _{x=1}=\Psi \left( x_{\widehat{T}},y_{\widehat{T}%
}\right) |_{_{\Gamma _{0},_{T}}},\text{ }\partial _{\nu }u\mid _{\Gamma
_{1,T}}=0.  \label{4.6}
\end{equation}%
Note that the function $\Psi \left( x_{\widehat{T}},y_{\widehat{T}}\right) $
in (\ref{4.6}) is taken from (\ref{4.1}). The boundary conditions for $m$
are:%
\begin{equation}
\left. 
\begin{array}{c}
\left( \sigma _{1}^{2}(x,y)m\right) _{x}\mid _{x=0}=0, \\ 
\left( \sigma _{2}^{2}(x,y)m\right) _{y}\mid _{y=0}=0, \\ 
\left( \sigma _{2}^{2}(x,y)m\right) _{y}\mid _{y=1}=0, \\ 
m\mid _{_{\Gamma _{0},_{T}}}=m\mid _{x=1}=0.%
\end{array}%
\right.  \label{4.7}
\end{equation}%
If 
\begin{equation}
\left. 
\begin{array}{c}
\sigma _{1}^{2}(\mathbf{x})=\sigma _{1},\text{ }\sigma _{2}^{2}(\mathbf{x}%
)=\sigma _{2}\text{ near }\partial \Omega , \\ 
\text{ where }\sigma _{1}>0\text{,}\sigma _{2}>0\text{ are some numbers,}%
\end{array}%
\right.  \label{4.8}
\end{equation}%
then conditions (\ref{4.7}) can be simplified as%
\begin{equation}
m\mid _{_{\Gamma _{0},_{T}}}=m\mid _{x=1}=0,\text{ }\partial _{\nu }m\mid
_{_{\Gamma _{1},_{T}}}=0.  \label{4.9}
\end{equation}

In addition to (\ref{2.4}), (\ref{4.6}) and (\ref{4.7}), we impose below the
following conditions on some functions involved in equations (\ref{4.4}), (%
\ref{4.5})%
\begin{equation}
c\in C\left( \overline{Q}_{T}\right) ,  \label{4.09}
\end{equation}%
\begin{equation}
\sigma _{1}^{2}(\mathbf{x})\geq 2\sigma _{0},\text{ }\sigma _{2}^{2}(\mathbf{%
x})\geq 2\sigma _{0}\text{ in }\overline{\Omega },  \label{4.90}
\end{equation}%
\begin{equation}
\sigma _{1}^{2}(\mathbf{x}),\text{ }\sigma _{2}^{2}(\mathbf{x})\in
C^{1}\left( \overline{\Omega }\right) ,  \label{4.901}
\end{equation}%
\begin{equation}
\left\Vert \sigma _{1}^{2}\right\Vert _{C^{1}\left( \overline{\Omega }%
\right) },\left\Vert \sigma _{2}^{2}\right\Vert _{C^{1}\left( \overline{%
\Omega }\right) }\leq D,  \label{4.902}
\end{equation}%
\begin{equation}
\text{ }\frac{\varphi _{1}^{2}(\mathbf{x})}{a_{0}(\mathbf{x})},\text{ }\frac{%
\varphi _{2}^{2}(\mathbf{x})}{b_{0}(\mathbf{x})}\in C^{1}\left( \overline{%
\Omega }\right) ,  \label{4.91}
\end{equation}%
\begin{equation}
\left\Vert c\right\Vert _{C\left( \overline{Q}_{T}\right) },\left\Vert
g\right\Vert _{C^{1}\left( \mathbb{R}^{2}\right) },\left\Vert \frac{\varphi
_{1}^{2}}{a_{0}}\right\Vert _{_{C^{1}\left( \overline{\Omega }\right)
}},\left\Vert \frac{\varphi _{2}^{2}}{b_{0}}\right\Vert _{_{C^{1}\left( 
\overline{\Omega }\right) }}\leq D,  \label{4.92}
\end{equation}%
where $\sigma _{0}>0$ and $D>0$ are certain numbers.

In the retrospective problem, functions $u\left( \mathbf{x},T\right) $ and $%
m\left( \mathbf{x},T\right) $ are supposed to be given. And one wants to use
this information to find functions $u\left( \mathbf{x},t\right) $ and $%
m\left( \mathbf{x},t\right) $ for $\left( \mathbf{x},t\right) \in Q_{T}.$
Thus, functions $u\left( \mathbf{x},T\right) $ and $m\left( \mathbf{x}%
,T\right) $ are the input data here. However, a valuable question to address
is about the stability of the problem of the determination of functions $%
u\left( \mathbf{x},t\right) $ and $m\left( \mathbf{x},t\right) $ for $\left( 
\mathbf{x},t\right) \in Q_{T}$ with respect to the noise in the input data.
This is exactly the question we address below. We are ready now to state the
retrospective problem which we address in this paper.

\textbf{Retrospective Problem.} \emph{Assume that conditions (\ref{2.2})-(%
\ref{2.5}), (\ref{2.11}), (\ref{2.12}), (\ref{4.04}), (\ref{4.8}) and (\ref%
{4.09})-(\ref{4.92}) hold. Suppose that we have two pairs of functions }%
\[
\left( u_{1},m_{1}\right) ,\left( u_{2},m_{2}\right) \in C^{2,1}\left( 
\overline{Q}_{T}\right) 
\]%
\emph{satisfying equations (\ref{4.4}), (\ref{4.5}) and boundary conditions (%
\ref{4.6}), (\ref{4.9}). Let} 
\begin{equation}
\left. 
\begin{array}{c}
u_{1}\left( \mathbf{x},T\right) =u_{1T}\left( \mathbf{x}\right) ,\text{ }%
u_{2}\left( \mathbf{x},T\right) =u_{2T}\left( \mathbf{x}\right) , \\ 
m_{1}\left( \mathbf{x},T\right) =m_{1T}\left( \mathbf{x}\right) ,\text{ }%
m_{2}\left( \mathbf{x},T\right) =m_{2T}\left( \mathbf{x}\right) .%
\end{array}%
\right.  \label{4.10}
\end{equation}%
\emph{Denote }%
\begin{equation}
\widetilde{u}_{T}\left( \mathbf{x}\right) =u_{1T}\left( \mathbf{x}\right)
-u_{2T}\left( \mathbf{x}\right) ,\text{ }\widetilde{m}_{T}\left( \mathbf{x}%
\right) =m_{1T}\left( \mathbf{x}\right) -m_{2T}\left( \mathbf{x}\right) .
\label{4.12}
\end{equation}%
\emph{Let a sufficiently small number }$\delta \in \left( 0,1\right) $\emph{%
\ be the level of the error in the input data (\ref{4.10}), i.e. let}%
\begin{equation}
\left\Vert \widetilde{u}_{T}\right\Vert _{H^{1}\left( \Omega \right) }\leq
\delta ,  \label{4.13}
\end{equation}%
\begin{equation}
\left\Vert \widetilde{m}_{T}\right\Vert _{H^{1}\left( \Omega \right) }\leq
\delta .  \label{4.14}
\end{equation}%
\emph{Estimate certain norms of differences }$\widetilde{u},\widetilde{m},$%
\begin{equation}
\widetilde{u}=u_{1}-u_{2},\text{ }\widetilde{m}=m_{1}-m_{2}  \label{4.11}
\end{equation}%
\emph{\ via the number }$\delta $ in \emph{(\ref{4.13}), (\ref{4.14}). }

\section{Carleman Estimates}

\label{sec:5}

The first step in addressing the Retrospective Problem is to prove new
Carleman estimates for two parabolic operators with variable coefficients of
their principal parts. The assumption of variable coefficients is necessary
since functions $\sigma _{1}^{2}(\mathbf{x}),\sigma _{2}^{2}(\mathbf{x})$ in
equations (\ref{4.4}), (\ref{4.5}) are not constants. We now derive two
Carleman estimates in the $n-$D case, $n\geq 1$. Below in sections 5,6 $%
x=\left( x_{1},x_{2},...,x_{n}\right) \in \mathbb{R}^{n}.$ Since in (\ref%
{2.2}) $\Omega =\left( 0,1\right) \times \left( 0,1\right) $ is a square,
then, to simplify the presentation and keeping in mind (\ref{2.11}), (\ref%
{2.12}) as well as the same notation for $\Omega $, we assume below that $%
\Omega $ is a cube, 
\begin{equation}
\left. 
\begin{array}{c}
\Omega =\left\{ x:0<x_{i}<1,i=1,...,n\right\} , \\ 
\partial \Omega =\Gamma _{0}\cup \Gamma _{1}, \\ 
\Gamma _{0}=\left\{ x_{1}=1,x_{i}\in \left( 0,1\right) ,i=2,...,n\right\} .%
\end{array}%
\right.  \label{5.1}
\end{equation}%
Since Carleman estimates are independent on lower terms of PDE operators 
\cite[Lemma 2.1.1]{KL}, then we work now only with principal parts of
elliptic operators. Consider two sets of functions $\left( a_{ij}\left(
x\right) \right) _{i,j=1}^{n}$ satisfying the following conditions:%
\begin{equation}
a_{ij}\left( x\right) \in C^{1}\left( \overline{\Omega }\right)
,a_{ij}\left( x\right) =a_{ji}\left( x\right) ;\text{ }i,j=1,...,n.
\label{5.2}
\end{equation}%
\begin{equation}
\mu _{1}\left\vert \xi \right\vert ^{2}\leq
\sum\limits_{i,j=1}^{n}a_{ij}\left( x\right) \xi _{i}\xi _{j}\leq \mu
_{2}\left\vert \xi \right\vert ^{2},\text{ }\forall \xi \in \mathbb{R}^{n},%
\text{ }\forall x\in \overline{\Omega },  \label{5.3}
\end{equation}%
where two numbers $\mu _{1},\mu _{2}>0$. Hence, we define two the elliptic
operator $L$ of the second order in the domain $\Omega $ as:%
\begin{equation}
Lu=\sum\limits_{i,j=1}^{n}\left( a_{ij}\left( x\right) u_{x_{j}}\right)
_{x_{i}}.  \label{5.6}
\end{equation}%
We define the normal derivative $\partial /\partial N$ at any side of the
cube (\ref{5.1}) as \cite[\S 1 of chapter 2]{Lad}%
\begin{equation}
\frac{\partial u}{\partial N}=\sum\limits_{j=1}^{n}a_{ij}\left( x\right)
u_{x_{j}}\cos \left( n,x_{i}\right) ,\text{ }i=1,..,n.  \label{5.8}
\end{equation}%
It follows from (\ref{5.8}) that 
\begin{equation}
\left. 
\begin{array}{c}
\text{if }Lu=\Delta u\text{ near }\partial \Omega , \\ 
\text{ then }\partial u/\partial N=\partial u/\partial \nu \text{ on }%
\partial \Omega .%
\end{array}%
\right.  \label{5.9}
\end{equation}%
We introduce the subspace $H_{0}^{2,1}\left( Q_{T}\right) $ of the space $%
H^{2,1}\left( Q_{T}\right) $ as:%
\begin{equation}
H_{0}^{2,1}\left( Q_{T}\right) =\left\{ u\in H^{2,1}\left( Q_{T}\right)
:u\mid _{\Gamma _{0,T}}=0,\frac{\partial u}{\partial N}\mid _{\Gamma
_{1,T}}=0\right\} .  \label{5.10}
\end{equation}

Let $\lambda >0$ and $s>0$ be two parameters, which we will choose later. We
introduce the Carleman Weight Function (CWF) $\varphi _{\lambda ,s}\left(
t\right) $ as \cite{MFG1}:%
\begin{equation}
\varphi _{\lambda ,s}\left( t\right) =e^{\lambda \left( t+2\right) ^{s}}.
\label{5.11}
\end{equation}

\subsection{Carleman estimate for the operator $\partial _{t}+L$}

\label{sec:5.1}

\textbf{Theorem 5.1}. \emph{Assume that conditions (\ref{5.2}), (\ref{5.3})
and (\ref{5.11}) hold. There exists a number }$C=C\left( T,\mu _{1}\right)
>0 $\emph{\ depending only on listed parameters, such that the following
Carleman estimate holds}%
\[
\int\limits_{Q_{T}}\left( u_{t}+Lu\right) ^{2}\varphi _{\lambda
,s}^{2}dxdt\geq \int\limits_{Q_{T}}\left( \frac{u_{t}^{2}}{4}+\left(
Lu\right) ^{2}\right) \varphi _{\lambda ,s}^{2}dxdt+ 
\]%
\begin{equation}
+C\lambda s\int\limits_{Q_{T}}\left( \nabla u\right) ^{2}\left( t+2\right)
^{s-1}\varphi _{\lambda ,s}^{2}dxdt+\frac{1}{2}\lambda
^{2}s^{2}\int\limits_{Q_{T}}\left( t+2\right) ^{2s-2}u^{2}\varphi _{\lambda
,s}^{2}-  \label{5.12}
\end{equation}%
\[
-\lambda s\left( T+2\right) ^{s-1}e^{2\lambda \left( T+2\right)
^{s}}\int\limits_{\Omega }u^{2}\left( x,T\right) dx-\mu _{1}e^{2\lambda
\left( T+2\right) ^{2}}\int\limits_{\Omega }\left( \nabla u\left( x,T\right)
\right) ^{2}dx, 
\]%
\[
\forall u\in H_{0}^{2,1}\left( Q_{T}\right) ,\text{ }\forall \lambda
>0,\forall s>1. 
\]

\textbf{Proof}. Below in section 5 $C=C\left( T,\mu _{1}\right) >0$\emph{\ }%
denotes different numbers depending only on listed parameters. Denote 
\begin{equation}
v=u\varphi _{\lambda ,s}=ue^{\lambda \left( t+2\right) ^{s}}.  \label{5.13}
\end{equation}%
Then 
\begin{equation}
\left. 
\begin{array}{c}
u=ve^{-\lambda \left( t+2\right) ^{s}}, \\ 
u_{t}=\left( v_{t}-\lambda s\left( t+2\right) ^{s-1}v\right) e^{-\lambda
\left( t+2\right) ^{s}}, \\ 
Lu=\left( Lv\right) e^{-\lambda \left( t+2\right) ^{s}}.%
\end{array}%
\right.  \label{5.14}
\end{equation}%
By (\ref{5.13}) and (\ref{5.14})%
\begin{equation}
\left. 
\begin{array}{c}
\left( u_{t}+Lu\right) ^{2}\varphi _{\lambda ,s}^{2}=\left[ \lambda s\left(
t+2\right) ^{s-1}v-\left( v_{t}+Lv\right) \right] ^{2}= \\ 
=\lambda ^{2}s^{2}\left( t+2\right) ^{2s-2}v^{2}-2\lambda s\left( t+2\right)
^{s-1}v\left( v_{t}+Lv\right) + \\ 
+\left( v_{t}+Lv\right) ^{2}.%
\end{array}%
\right.  \label{5.15}
\end{equation}

\emph{Step 1.} Estimate from the below the term $-2\lambda s\left(
t+2\right) ^{s-1}v\left( v_{t}+Lv\right) $ in (\ref{5.15}),%
\[
-2\lambda s\left( t+2\right) ^{s-1}v\left( v_{t}+Lv\right) =\left( -\lambda
s\left( t+2\right) ^{s-1}v^{2}\right) _{t}+\lambda s\left( s-1\right) \left(
t+2\right) ^{s-2}v^{2}- 
\]%
\begin{equation}
-2\lambda s\left( t+2\right) ^{s-1}\sum\limits_{i,j=1}^{n}\left(
a_{ij}\left( x\right) v_{x_{j}}\right) _{x_{i}}v\geq  \label{5.16}
\end{equation}%
\[
\geq \left( -\lambda s\left( t+2\right) ^{s-1}v^{2}\right) _{t}-2\lambda
s\left( t+2\right) ^{s-1}\sum\limits_{i,j=1}^{n}\left( a_{ij}\left( x\right)
v_{x_{j}}\right) _{x_{i}}v. 
\]%
Using (\ref{5.3}), estimate now the second term in the third line of (\ref%
{5.16}),%
\[
-2\lambda s\left( t+2\right) ^{s-1}\sum\limits_{i,j=1}^{n}\left(
a_{ij}\left( x\right) v_{x_{j}}\right)
_{x_{i}}v=\sum\limits_{i,j=1}^{n}\left( -2\lambda s\left( t+2\right)
^{s-1}a_{ij}\left( x\right) v_{x_{j}}v\right) _{x_{i}}+ 
\]%
\[
+2\lambda s\left( t+2\right) ^{s-1}\sum\limits_{i,j=1}^{n}a_{ij}\left(
x\right) v_{x_{j}}v_{x_{i}}\geq 
\]%
\[
\geq \sum\limits_{i,j=1}^{n}\left( -2\lambda s\left( t+2\right)
^{s-1}a_{ij}\left( x\right) v_{x_{j}}v\right) _{x_{i}}+2\mu _{1}\lambda
s\left( t+2\right) ^{s-1}\left\vert \nabla v\right\vert ^{2}. 
\]%
Comparing this with (\ref{5.15}) and (\ref{5.16}), we obtain%
\[
\left( u_{t}+L_{1}u\right) ^{2}\varphi _{\lambda ,s}^{2}\geq \lambda
^{2}s^{2}\left( t+2\right) ^{2s-2}v^{2}+C\lambda s\left( t+2\right)
^{s-1}\left\vert \nabla v\right\vert ^{2}+ 
\]%
\begin{equation}
+\left( v_{t}+Lv\right) ^{2}+  \label{5.17}
\end{equation}%
\[
+\left( -\lambda s\left( t+2\right) ^{s-1}v^{2}\right)
_{t}+\sum\limits_{i,j=1}^{n}\left( -2\lambda s\left( t+2\right)
^{s-1}a_{ij}\left( x\right) v_{x_{j}}v\right) _{x_{i}}. 
\]

\emph{Step 2}. Evaluate the term $\left( v_{t}+Lv\right) ^{2}$ in (\ref{5.17}%
). Using (\ref{5.6}), we obtain%
\[
\left( v_{t}+Lv\right) ^{2}=v_{t}^{2}+2v_{t}Lv+\left( Lv\right) ^{2}= 
\]%
\begin{equation}
+v_{t}^{2}+\left( Lv\right) ^{2}+2\sum\limits_{i,j=1}^{n}\left( a_{ij}\left(
x\right) v_{x_{j}}\right) _{x_{i}}v_{t}.  \label{5.18}
\end{equation}%
Estimate the third time in the second line of (\ref{5.18}), 
\begin{equation}
2\sum\limits_{i,j=1}^{n}\left( a_{ij}\left( x\right) v_{x_{j}}\right)
_{x_{i}}v_{t}=\sum\limits_{i,j=1}^{n}\left( 2a_{ij}\left( x\right)
v_{x_{j}}v_{t}\right) _{x_{i}}-2\sum\limits_{i,j=1}^{n}a_{ij}\left( x\right)
v_{x_{j}}v_{tx_{i}}.  \label{5.19}
\end{equation}%
Since by (\ref{5.2}) $a_{ij}\left( x\right) =a_{ji}\left( x\right) ,$ then
the last term of (\ref{5.19}) is: 
\[
-2\sum\limits_{i,j=1}^{n}a_{ij}\left( x\right)
v_{x_{j}}v_{tx_{i}}=-\sum\limits_{i,j=1}^{n}a_{ij}\left( x\right) \left(
v_{x_{j}}v_{tx_{i}}+v_{x_{j}}v_{tx_{i}}\right) = 
\]%
\begin{equation}
=\left( -\sum\limits_{i,j=1}^{n}a_{ij}\left( x\right)
v_{x_{j}}v_{x_{i}}\right) _{t}.  \label{5.20}
\end{equation}%
Thus, we have proven in (\ref{5.19}) and (\ref{5.20}) that 
\begin{equation}
2\sum\limits_{i,j=1}^{n}\left( a_{ij}\left( x\right) v_{x_{j}}\right)
_{x_{i}}v_{t}=\sum\limits_{i,j=1}^{n}\left( 2a_{ij}\left( x\right)
v_{x_{j}}v_{t}\right) _{x_{i}}+\left( -\sum\limits_{i,j=1}^{n}a_{ij}\left(
x\right) v_{x_{j}}v_{x_{i}}\right) _{t}.  \label{5.200}
\end{equation}%
Thus, (\ref{5.18})-(\ref{5.200}) imply%
\begin{equation}
\left( v_{t}+Lv\right) ^{2}=v_{t}^{2}+\left( Lv\right)
^{2}+\sum\limits_{i,j=1}^{n}\left( 2a_{ij}\left( x\right)
v_{x_{j}}v_{t}\right) _{x_{i}}+  \label{5.21}
\end{equation}%
\[
+\sum\limits_{i,j=1}^{n}\left( 2a_{ij}\left( x\right) v_{x_{j}}v_{t}\right)
_{x_{i}}+\left( -\sum\limits_{i,j=1}^{n}a_{ij}\left( x\right)
v_{x_{j}}v_{x_{i}}\right) _{t}. 
\]%
Therefore, combining (\ref{5.21}) with (\ref{5.17}), we obtain%
\[
\left( u_{t}+Lu\right) ^{2}\varphi _{\lambda ,s}^{2}\geq v_{t}^{2}+\left(
Lv\right) ^{2}+\lambda ^{2}s^{2}\left( t+2\right) ^{2s-2}v^{2}+C\lambda
s\left( t+2\right) ^{s-1}\left\vert \nabla v\right\vert ^{2}+ 
\]%
\begin{equation}
+\sum\limits_{i,j=1}^{n}\left( -2\lambda s\left( t+2\right)
^{s-1}a_{ij}\left( x\right) v_{x_{j}}v\right) _{x_{i}}+\left(
2\sum\limits_{i,j=1}^{n}a_{ij}\left( x\right) v_{x_{j}}v_{t}\right) _{x_{i}}
\label{5.22}
\end{equation}%
\[
+\left( -\lambda s\left( t+2\right)
^{s-1}v^{2}-\sum\limits_{i,j=1}^{n}a_{ij}\left( x\right)
v_{x_{j}}v_{x_{i}}\right) _{t}. 
\]

Using (\ref{5.13}), we now need to replace the function $v$ in (\ref{5.22})
with the function $u.$ We want to keep terms with $u_{t}^{2}$ and $u^{2}$
with positive signs at them. Hence, we put a special attention to the term $%
v_{t}^{2}+\lambda ^{2}s^{2}\left( t+2\right) ^{2s-2}v^{2}$ in the first line
of (\ref{5.22}). Let $c\in \left( 0,1\right) $ be a number, which we will
choose later. Using (\ref{5.13}) and Cauchy-Schwarz inequality, we obtain%
\[
v_{t}^{2}+\lambda ^{2}s^{2}\left( t+a\right) ^{2s-2}v^{2}\geq
cv_{t}^{2}+\lambda ^{2}s^{2}\left( t+2\right) ^{2s-2}v^{2}= 
\]%
\[
=c\left( u_{t}+\lambda s\left( t+2\right) ^{s-1}u\right) ^{2}e^{2\lambda
\left( t+2\right) ^{s}}= 
\]%
\[
=c\left( u_{t}^{2}+2u_{t}\cdot \lambda s\left( t+2\right) ^{s-1}u+\lambda
^{2}s^{2}\left( t+2\right) ^{2s-2}u^{2}\right) e^{2\lambda \left( t+2\right)
^{s}}+ 
\]%
\[
+\lambda ^{2}s^{2}\left( t+2\right) ^{2s-2}u^{2}e^{2\lambda \left(
t+2\right) ^{s}}\geq 
\]%
\[
\geq c\left( \frac{u_{t}^{2}}{2}-\lambda ^{2}s^{2}\left( t+2\right)
^{2s-2}u^{2}\right) e^{2\lambda \left( t+2\right) ^{s}}+\lambda
^{2}s^{2}\left( t+2\right) ^{2s-2}u^{2}e^{2\lambda \left( t+2\right) ^{s}}= 
\]%
\[
=\frac{c}{2}u_{t}^{2}e^{2\lambda \left( t+2\right) ^{s}}+\left( 1-c\right)
\lambda ^{2}s^{2}\left( t+2\right) ^{2s-2}u^{2}e^{2\lambda \left( t+2\right)
^{s}}. 
\]%
Choosing $c=1/2,$ we obtain%
\[
v_{t}^{2}+\lambda ^{2}s^{2}\left( t+2\right) ^{2s-2}v^{2}\geq \frac{1}{4}%
u_{t}^{2}e^{2\lambda \left( t+2\right) ^{s}}+\frac{1}{2}\lambda
^{2}s^{2}\left( t+2\right) ^{2s-2}u^{2}e^{2\lambda \left( t+2\right) ^{s}}. 
\]%
Combining this with (\ref{5.22}), we obtain the following pointwise Carleman
estimate for the operator $u_{t}+Lu:$%
\[
\left( u_{t}+Lu\right) ^{2}\varphi _{\lambda ,s}^{2}\geq \left( \frac{1}{4}%
u_{t}^{2}+\left( Lu\right) ^{2}\right) \varphi _{\lambda ,s}^{2}+ 
\]%
\[
+C\lambda s\left( t+2\right) ^{s-1}\left\vert \nabla u\right\vert
^{2}\varphi _{\lambda ,s}^{2}+\frac{1}{2}\lambda ^{2}s^{2}\left( t+2\right)
^{2s-2}u^{2}\varphi _{\lambda ,s}^{2}+ 
\]%
\[
+\sum\limits_{i,j=1}^{n}\left( -2\lambda s\left( t+2\right)
^{s-1}a_{ij}\left( x\right) u_{x_{j}}u\varphi _{\lambda ,s}^{2}\right)
_{x_{i}}+ 
\]%
\begin{equation}
+\left( 2\sum\limits_{i,j=1}^{n}a_{ij}\left( x\right) u_{x_{j}}\left(
u_{t}+\lambda s\left( t+2\right) ^{s-1}\right) \varphi _{\lambda
,s}^{2}\right) _{x_{i}}+  \label{5.23}
\end{equation}%
\[
+\left( -\lambda s\left( t+2\right) ^{s-1}u^{2}\varphi _{\lambda
,s}^{2}-\sum\limits_{i,j=1}^{n}a_{ij}\left( x\right)
u_{x_{j}}u_{x_{i}}\varphi _{\lambda ,s}^{2}\right) _{t}. 
\]%
Integrate inequality (\ref{5.23}) over $Q_{T}$. Using Gauss formula and (\ref%
{5.8})-(\ref{5.10}), we obtain (\ref{5.12}), which is the target estimate of
this theorem. $\square $

\subsection{Carleman estimate for the operator $\partial _{t}-L$}

\label{sec:5.2}

\textbf{Theorem 5.2.} \emph{Assume that conditions of Theorem 5.1 hold. Then
there exists a number }$C=C\left( T,\mu _{1}\right) >0$\emph{\ depending
only on listed parameters and a sufficiently large absolute number }$s_{0}>1$%
\emph{\ such that the following Carleman estimate holds:}%
\[
\int\limits_{Q_{T}}\left( u_{t}-Lu\right) ^{2}\varphi _{\lambda
,s}^{2}dxdt\geq 
\]%
\begin{equation}
\geq \mu _{1}\sqrt{s}\int\limits_{Q_{Y}}\left( \nabla u\right) ^{2}\varphi
_{\lambda ,s}^{2}dxdt+C\lambda s^{2}\int\limits_{Q_{T}}\left( t+2\right)
^{s-1}u^{2}\varphi _{\lambda ,s}^{2}dxdt-  \label{5.24}
\end{equation}%
\[
-\lambda s\left( T+2\right) ^{s-1}e^{2\lambda \left( T+2\right)
^{s}}\int\limits_{\Omega }u^{2}\left( x,T\right) dx- 
\]%
\[
-e^{2^{s+1}\lambda }\int\limits_{\Omega }\left( \mu _{2}\left( \nabla
u\left( x,0\right) \right) ^{2}+\frac{\sqrt{s}}{2}u^{2}\left( x,0\right)
\right) dx, 
\]%
\[
\forall u\in H_{0}^{2,1}\left( Q_{T}\right) ,\text{ }\forall \lambda
>0,\forall s>s_{0}. 
\]

\textbf{Proof}. Just as in the proof of Theorem 5.1, introduce the new
function $v$ as in (\ref{5.13}). Hence, using (\ref{5.14}), we obtain%
\[
\left( u_{t}-Lu\right) ^{2}\varphi _{\lambda ,s}^{2}=\left[ v_{t}-\left(
\lambda s\left( t+2\right) ^{s-1}v+Lv\right) \right] ^{2}\geq 
\]%
\[
\geq -2v_{t}\left( \lambda s\left( t+2\right) ^{s-1}v+Lv\right) = 
\]%
\begin{equation}
=\left( -\lambda s\left( t+2\right) ^{s-1}v^{2}\right) _{t}+\lambda s\left(
s-1\right) \left( t+2\right) ^{s-1}v^{2}-  \label{5.25}
\end{equation}%
\[
-2\sum\limits_{i,j=1}^{n}\left( a_{ij}\left( x\right) v_{x_{j}}\right)
_{x_{i}}v_{t}. 
\]%
By the last line of (\ref{5.25}) is%
\begin{equation}
-2\sum\limits_{i,j=1}^{n}\left( a_{ij}\left( x\right) v_{x_{j}}\right)
_{x_{i}}v_{t}=\sum\limits_{i,j=1}^{n}\left( -2a_{ij}\left( x\right)
v_{x_{j}}v_{t}\right) _{x_{i}}+  \label{5.26}
\end{equation}%
\[
+\left( \sum\limits_{i,j=1}^{n}a_{ij}\left( x\right)
v_{x_{j}}v_{x_{i}}\right) _{t}. 
\]%
Hence, using (\ref{5.25}) and (\ref{5.26}), we obtain 
\[
\left( u_{t}-Lu\right) ^{2}\varphi _{\lambda ,s}^{2}\geq \lambda s\left(
s-1\right) \left( t+2\right) ^{s-1}u^{2}\varphi _{\lambda ,s}^{2}+ 
\]%
\[
+\left( -\lambda s\left( t+2\right) ^{s-1}u^{2}\varphi _{\lambda
,s}^{2}+\sum\limits_{i,j=1}^{n}a_{ij}\left( x\right)
u_{x_{j}}u_{x_{i}}\varphi _{\lambda ,s}^{2}\right) _{t}+ 
\]%
\[
+\sum\limits_{i,j=1}^{n}\left( -2a_{ij}\left( x\right) u_{x_{j}}\left(
u_{t}+\lambda s\left( t+2\right) ^{s-1}u\right) \varphi _{\lambda
,s}^{2}\right) _{x_{i}}. 
\]%
Integrate this inequality over $Q_{T}$ and use Gauss formula as well as (\ref%
{5.10}). We obtain%
\begin{equation}
\int\limits_{Q_{T}}\left( u_{t}-Lu\right) ^{2}\varphi _{\lambda
,s}^{2}dxdt\geq \lambda s\left( s-1\right) \int\limits_{Q_{T}}\left(
t+2\right) ^{s-1}u^{2}\varphi _{\lambda ,s}^{2}dxdt-  \label{5.27}
\end{equation}%
\[
-\lambda s\left( T+2\right) ^{s-1}e^{2\lambda \left( T+2\right)
^{s}}\int\limits_{\Omega }u^{2}\left( x,T\right) dx-\mu
_{2}e^{2^{s+1}\lambda }\int\limits_{\Omega }\left( \nabla u\left( x,0\right)
\right) ^{2}dx. 
\]

We now need to incorporate in the first line of (\ref{5.27}) the term with $%
\left( \nabla u\right) ^{2}$ with the positive sign at it. To do this,
consider the following expression:%
\[
\left( u_{t}-Lu\right) u\varphi _{\lambda ,s}^{2}=\left( \frac{u^{2}}{2}%
\varphi _{\lambda ,s}^{2}\right) _{t}-\lambda s\left( t+2\right)
^{s-1}u^{2}\varphi _{\lambda ,s}^{2}- 
\]%
\[
-\sum\limits_{i,j=1}^{n}\left( a_{ij}\left( x\right) u_{x_{j}}\right)
_{x_{i}}u\varphi _{\lambda ,s}^{2}= 
\]%
\[
=\left( \frac{u^{2}}{2}\varphi _{\lambda ,s}^{2}\right) _{t}-\lambda s\left(
t+2\right) ^{s-1}u^{2}\varphi _{\lambda
,s}^{2}+\sum\limits_{i,j=1}^{n}\left( -a_{ij}\left( x\right)
u_{x_{j}}u\varphi _{\lambda ,s}^{2}\right) _{x_{i}}+ 
\]%
\[
+\sum\limits_{i,j=1}^{n}a_{ij}\left( x\right) u_{x_{j}}u_{x_{i}}\varphi
_{\lambda ,s}^{2}. 
\]%
Thus, 
\[
\left( u_{t}-Lu\right) u\varphi _{\lambda ,s}^{2}\geq \mu _{1}\left( \nabla
u\right) ^{2}\varphi _{\lambda ,s}^{2}-\lambda s\left( t+a\right)
^{s-1}u^{2}\varphi _{\lambda ,s}^{2}+ 
\]%
\[
+\left( \frac{u^{2}}{2}\varphi _{\lambda ,s}^{2}\right)
_{t}+\sum\limits_{i,j=1}^{n}\left( -a_{ij}\left( x\right) u_{x_{j}}u\varphi
_{\lambda ,s}^{2}\right) _{x_{i}}. 
\]%
Integrating this inequality over $Q_{T}$ and using Gauss formula and (\ref%
{5.10}), we obtain%
\[
\int\limits_{Q_{Y}}\left( u_{t}-Lu\right) u\varphi _{\lambda ,s}^{2}dxdt\geq
\mu _{1}\int\limits_{Q_{Y}}\left( \nabla u\right) ^{2}\varphi _{\lambda
,s}^{2}dxdt-\lambda s\int\limits_{Q_{Y}}\left( t+2\right) ^{s-1}u^{2}\varphi
_{\lambda ,s}^{2}dxdt- 
\]%
\[
-\frac{1}{2}e^{2^{s+1}\lambda }\int\limits_{\Omega }u^{2}\left( x,0\right)
dx. 
\]%
Multiply this inequality by $\sqrt{s}>\sqrt{s_{0}}$ and sum up with (\ref%
{5.27}). We obtain%
\[
\sqrt{s}\int\limits_{Q_{Y}}\left( u_{t}-Lu\right) u\varphi _{\lambda
,s}^{2}dxdt+\int\limits_{Q_{T}}\left( u_{t}-Lu\right) ^{2}\varphi _{\lambda
,s}^{2}dxdt\geq 
\]%
\begin{equation}
\geq \mu _{1}\sqrt{s}\int\limits_{Q_{Y}}\left( \nabla u\right) ^{2}\varphi
_{\lambda ,s}^{2}dxdt+C\lambda s^{2}\int\limits_{Q_{T}}\left( t+2\right)
^{s-1}u^{2}\varphi _{\lambda ,s}^{2}dxdt-  \label{5.28}
\end{equation}%
\[
-\lambda s\left( T+2\right) ^{s-1}e^{2\lambda \left( T+2\right)
^{s}}\int\limits_{\Omega }u^{2}\left( x,T\right) dx- 
\]%
\[
-e^{2^{s+1}\lambda }\int\limits_{\Omega }\left( \mu _{2}\left( \nabla
u\left( x,0\right) \right) ^{2}+\frac{\sqrt{s}}{2}u^{2}\left( x,0\right)
\right) dx. 
\]%
Finally, noting that 
\[
\sqrt{s}\int\limits_{Q_{Y}}\left( u_{t}-Lu\right) u\varphi _{\lambda
,s}^{2}dxdt+\int\limits_{Q_{T}}\left( u_{t}-Lu\right) ^{2}\varphi _{\lambda
,s}^{2}dxdt\leq 
\]%
\[
\leq 2\int\limits_{Q_{T}}\left( u_{t}-Lu\right) ^{2}\varphi _{\lambda
,s}^{2}dxdt+\frac{s}{2}\int\limits_{Q_{T}}u^{2}\varphi _{\lambda ,s}^{2}dxdt 
\]%
and substituting this in (\ref{5.28}), we obtain (\ref{5.24}), which is the
target estimate of this theorem. $\square $

\section{A Generalized Retrospective Problem}

\label{sec:6}

In this section we obtain H\"{o}lder stability estimate for a generalized
retrospective problem for a generalized MFGS. This will help us in the next
section to obtain H\"{o}lder stability estimate for our target Retrospective
Problem posed in section 4. Notations of section 5 are kept here.

Let functions 
\begin{equation}
F_{1}\in C^{1}\left( \mathbb{R}^{3n+4}\right) ,\text{ }F_{2}\in C^{1}\left( 
\mathbb{R}^{2n+5}\right) .  \label{6.1}
\end{equation}%
Also, let functions

\begin{equation}
K_{1}(x),\dots ,K_{n}(x)\in L_{\infty }(\Omega ).  \label{6.2}
\end{equation}

We consider the form of the MFGS, which is slightly more general than the
conventional one of, e.g. \cite{A}. The first equation of this system is:

\begin{equation}
\left. 
\begin{array}{c}
u_{t}+Lu+ \\ 
+F_{1}(u,\nabla u,m,\nabla m,\int\limits_{0}^{1}K_{1}(x)m(x,t)dx_{1}, \\ 
\int\limits_{0}^{1}K_{2}(x)m(x,t)dx_{2}...,,\int%
\limits_{0}^{1}K_{n}(x)m(x,t)dx_{n},x,t)=0%
\end{array}%
\right.  \label{6.3}
\end{equation}



The second equation is:%
\begin{equation}
m_{t}-Lm+F_{2}\left( \nabla m,m,u,\nabla u,Lu,x,t\right) =0.  \label{6.4}
\end{equation}%
The boundary conditions now are similar with ones of (\ref{4.6}), (\ref{4.9}%
),%
\begin{equation}
\left. 
\begin{array}{c}
u\mid _{\Gamma _{0,T}}=f_{u}\left( x,t\right) ,\text{ }\partial u/\partial
N\mid _{\Gamma _{1,T}}=0, \\ 
m\mid _{\Gamma _{0,T}}=f_{m}\left( x,t\right) ,\text{ }\partial m/\partial
N\mid _{\Gamma _{1,T}}=0,%
\end{array}%
\right.  \label{6.5}
\end{equation}%
where $f_{u}$ and $f_{m}$ are certain functions defined on $\Gamma _{0,T}.$

\textbf{Generalized Retrospective Problem.} \emph{Assume that conditions (%
\ref{6.1}) and (\ref{6.2}) hold. Suppose that we have two pairs of functions 
}%
\[
\left( u_{1},m_{1}\right) ,\left( u_{2},m_{2}\right) \in C^{2,1}\left( 
\overline{Q}_{T}\right) 
\]%
\emph{satisfying (\ref{6.3})-(\ref{6.5}). Let} 
\begin{equation}
\left. 
\begin{array}{c}
u_{1}\left( x,T\right) =u_{1T}\left( x\right) ,u_{2}\left( x,T\right)
=u_{2T}\left( x\right) , \\ 
m_{1}\left( x,T\right) =m_{1T}\left( x\right) ,m_{2}\left( x,T\right)
=m_{2T}\left( x\right) .%
\end{array}%
\right.  \label{6.6}
\end{equation}%
Following \emph{(\ref{4.12})-(\ref{4.11}),} \emph{denote }%
\begin{equation}
\widetilde{u}_{T}\left( x\right) =u_{1T}\left( x\right) -u_{2T}\left(
x\right) ,\text{ }\widetilde{m}_{T}\left( x\right) =m_{1T}\left( x\right)
-m_{2T}\left( x\right) .  \label{6.06}
\end{equation}%
\emph{Let a sufficiently small number }$\delta \in \left( 0,1\right) $\emph{%
\ be the level of the error in the input data (\ref{6.6}), i.e. let}%
\begin{equation}
\left\Vert \widetilde{u}_{T}\right\Vert _{H^{1}\left( \Omega \right) }\leq
\delta ,  \label{6.60}
\end{equation}%
\begin{equation}
\left\Vert \widetilde{m}_{T}\right\Vert _{H^{1}\left( \Omega \right) }\leq
\delta .  \label{6.61}
\end{equation}%
\emph{Estimate certain norms of differences }$\widetilde{u},\widetilde{m},$%
\begin{equation}
\widetilde{u}=u_{1}-u_{2},\text{ }\widetilde{m}=m_{1}-m_{2}  \label{6.7}
\end{equation}%
\emph{\ via the number }$\delta $ in \emph{(\ref{6.60}), (\ref{6.61}).}

Let the number $M_{1}>0$. Consider the set $B\left( M_{1}\right) $ of pairs
of functions $\left( u,m\right) $ such that%
\begin{equation}
B\left( M_{1}\right) =\left\{ 
\begin{array}{c}
\left( u,m\right) \in C^{2,1}\left( \overline{Q}_{T}\right) : \\ 
\max_{t\in \left[ 0,T\right] }\left\Vert u\left( x,t\right) \right\Vert
_{C^{2}\left( \overline{\Omega }\right) }\leq M_{1}, \\ 
\max_{t\in \left[ 0,T\right] }\left\Vert m\left( x,t\right) \right\Vert
_{C^{1}\left( \overline{\Omega }\right) }\leq M_{1}%
\end{array}%
\right\} .  \label{6.9}
\end{equation}%
Also, let 
\begin{equation}
\left\Vert K_{1}\right\Vert _{L_{\infty }\left( \Omega \right) },\dots
,\left\Vert K_{n}\right\Vert _{L_{\infty }\left( \Omega \right) }\leq M_{1}.
\label{6.10}
\end{equation}%
It follows from (\ref{6.1}), (\ref{6.2}) and (\ref{6.9}) that there exists a
number

$M=M\left( M_{1},K_{1},K_{2},...,K_{n},Q_{T}\right) >0$ such that%
\begin{equation}
\left. 
\begin{array}{c}
\max_{\left( u,m,x,t\right) \in \overline{B\left( M_{1}\right) }\times 
\overline{Q}_{T}} \\ 
\left\vert 
\begin{array}{c}
\widehat{\nabla }F_{1}(u,\nabla u,m,\nabla
m,\int\limits_{0}^{1}K_{1}(x)m(x,t)dx_{1},\int%
\limits_{0}^{1}K_{2}(x)m(x,t)dx_{2},... \\ 
\dots ,\int\limits_{0}^{1}K_{n}(x)m(x,t)dx_{n},x,t)%
\end{array}%
\right\vert \leq \\ 
\leq M, \\ 
\max_{\left( u,m,x,t\right) \in \overline{B\left( M_{1}\right) }\times 
\overline{Q}_{T}}\left\vert \widehat{\nabla }F_{2}\left( u,\nabla u,m,\nabla
m,Lu,x,t\right) \right\vert \leq M,%
\end{array}%
\right.  \label{6.11}
\end{equation}%
where \textquotedblleft $\widehat{\nabla }"$ means that the first
derivatives are taken only with respect to those variables, which contain $%
u,\nabla u,m$ and $\nabla m.$ For any number $\gamma \in \left( 0,T\right) $
denote 
\begin{equation}
Q_{\gamma T}=\Omega \times \left( \gamma ,T\right) \subset Q_{T}.
\label{6.12}
\end{equation}

Let a sufficiently small number $\delta \in \left( 0,1\right) $
characterizes the level of the noise in the input data

\textbf{Theorem 6.1 }(H\"{o}lder stability estimate). \emph{Assume that
conditions (\ref{6.6})-(\ref{6.12}) hold and let two pairs of functions }$%
\left( u_{1},m_{1}\right) ,\left( u_{2},m_{2}\right) \in B\left(
M_{1}\right) $\emph{. satisfy equations (\ref{6.3}), (\ref{6.4}) and
boundary conditions (\ref{6.5}). Suppose that (\ref{6.60}) and (\ref{6.61})
hold. Let }$T>1$\emph{\ and in (\ref{6.12})} 
\begin{equation}
\gamma \in \left( 1,T\right) .  \label{6.122}
\end{equation}%
\emph{Then there exists a sufficiently small number }$\delta _{0}=\delta
_{0}\left( T,M_{1},M,\gamma \right) \in \left( 0,1\right) $\emph{\ and a
number }$\beta =\beta \left( T,M_{1},M,\gamma \right) \in \left(
0,1/6\right) ,$\emph{\ both these numbers depend only on listed parameters,
such that the following H\"{o}lder stability estimate for the Generalized
Retrospective Problem is valid:}%
\begin{equation}
\left. 
\begin{array}{c}
\left\Vert \widetilde{u}_{t}\right\Vert _{L_{2}\left( Q_{\gamma T}\right)
}+\left\Vert L\widetilde{u}\right\Vert _{L_{2}\left( Q_{\gamma T}\right)
}+\left\Vert \nabla \widetilde{u}\right\Vert _{L_{2}\left( Q_{\gamma
T}\right) }+\left\Vert \widetilde{u}\right\Vert _{L_{2}\left( Q_{\gamma
T}\right) }+ \\ 
+\left\Vert \nabla \widetilde{m}\right\Vert _{L_{2}\left( Q_{\gamma
T}\right) }+\left\Vert \widetilde{m}\right\Vert _{L_{2}\left( Q_{\gamma
T}\right) }\leq C_{1}\delta ^{\beta },\text{ }\forall \delta \in \left(
0,\delta _{0}\right) ,%
\end{array}%
\right.  \label{6.123}
\end{equation}%
\emph{where the number }$C_{1}=C_{1}\left( T,M_{1},M,\gamma \right) >0$\emph{%
\ depends only on listed parameters. In addition, problem (\ref{6.3})-(\ref%
{6.5}) has at most one solution }$\left( u,m\right) \in B\left( M_{1}\right)
.$

\textbf{Proof}. Below in section 6 $C_{1}=C_{1}\left( T,M_{1},M,\gamma
\right) >0$ denotes different numbers depending only on listed parameters.
First, write equations (\ref{6.3}) and (\ref{6.4}) for the pair $\left(
u_{1},m_{1}\right) .$ Then repeat this for the pair $\left(
u_{2},m_{2}\right) .$ Next, subtract equations for $\left(
u_{2},m_{2}\right) $ from equations for $\left( u_{1},m_{1}\right) .$ Using
the multidimensional analog of Taylor formula \cite{V} and (\ref{6.11}), we
obtain a system of two equations with respect to the pair $\left( \widetilde{%
u},\widetilde{m}\right) ,$ see (\ref{6.7}) for this pair. However, to
simplify the presentation, we turn this system in two differential
inequalities, as it is conventionally done when Carleman estimates are
applied, see, e.g. books \cite{KL,MFGbook,LRS}. This system is:%
\begin{equation}
\left\vert \widetilde{u}_{t}+L\widetilde{u}\right\vert \leq C_{1}\left(
\left\vert \nabla \widetilde{u}\right\vert +\left\vert \widetilde{u}%
\right\vert +\left\vert \nabla \widetilde{m}\right\vert +\left\vert 
\widetilde{m}\right\vert +\int\limits_{\Omega }\left\vert \widetilde{m}%
\left( y,t\right) \right\vert dy\right) \text{ in }Q_{T},  \label{6.13}
\end{equation}%
\begin{equation}
\left\vert \widetilde{m}_{t}-L\widetilde{m}\right\vert \leq C_{1}\left(
\left\vert \nabla \widetilde{m}\right\vert +\left\vert \widetilde{m}%
\right\vert +\left\vert \nabla \widetilde{u}\right\vert +\left\vert 
\widetilde{u}\right\vert +\left\vert L\widetilde{u}\right\vert \right) \text{
in }Q_{T}.  \label{6.14}
\end{equation}%
In addition, by (\ref{5.10}), (\ref{6.5}), (\ref{6.6}), (\ref{6.06}) and (%
\ref{6.7})%
\begin{equation}
\widetilde{u},\widetilde{m}\in H_{0}^{2,1}\left( Q_{T}\right) ,  \label{6.15}
\end{equation}%
\begin{equation}
\widetilde{u}\left( x,T\right) =\widetilde{u}_{T}\left( x\right) ,
\label{6.16}
\end{equation}%
\begin{equation}
\widetilde{m}\left( x,T\right) =\widetilde{m}_{T}\left( x\right) .
\label{6.17}
\end{equation}

First, we multiply both sides of inequality (\ref{6.13}) by the CWF (\ref%
{5.11}), square both sides of the resulting inequality, integrate over the
domain $Q_{T}$, use Carleman estimate (\ref{5.12}), (\ref{6.15}), (\ref{6.16}%
) and Cauchy-Schwarz inequality. We obtain%
\[
\int\limits_{Q_{T}}\left( \frac{\widetilde{u}_{t}^{2}}{4}+\left( L\widetilde{%
u}\right) ^{2}\right) \varphi _{\lambda ,s}^{2}dxdt+ 
\]%
\begin{equation}
+C\lambda s\int\limits_{Q_{T}}\left( \nabla \widetilde{u}\right) ^{2}\varphi
_{\lambda ,s}^{2}dxdt+\frac{1}{2}\lambda ^{2}s^{2}\int\limits_{Q_{T}}%
\widetilde{u}^{2}\varphi _{\lambda ,s}^{2}\leq  \label{6.18}
\end{equation}%
\[
\leq C_{1}\int\limits_{Q_{T}}\left( \left\vert \nabla \widetilde{u}%
\right\vert ^{2}+\widetilde{u}^{2}+\left\vert \nabla \widetilde{m}%
\right\vert ^{2}+\widetilde{m}^{2}\right) \varphi _{\lambda ,s}^{2}dxdt 
\]%
\[
+\lambda s\left( T+2\right) ^{s-1}e^{2\lambda \left( T+2\right)
^{s}}\int\limits_{\Omega }\widetilde{u}_{T}^{2}\left( x\right) dx+\mu
_{1}e^{2\lambda \left( T+2\right) ^{2}}\int\limits_{\Omega }\left( \nabla 
\widetilde{u}_{T}\left( x\right) \right) ^{2}dx, 
\]%
\[
\forall \lambda >0,\forall s>1. 
\]

We set from now on:%
\begin{equation}
\lambda \geq 1,s>s_{0},  \label{6.19}
\end{equation}%
where $s_{0}>>1$ is the number of Theorem 5.2. Choose the number $s_{1},$ 
\begin{equation}
s_{1}=s_{1}\left( T,M_{1},M,\gamma \right) \geq s_{0}  \label{6.20}
\end{equation}%
depending only on listed parameters and such that 
\begin{equation}
2C_{1}<\max \left( Cs,\frac{s^{2}}{2}\right) ,\text{ }\forall s\geq s_{1}.
\label{6.21}
\end{equation}%
Then, using (\ref{6.60}) and (\ref{6.18})-(\ref{6.21}), we obtain 
\[
\int\limits_{Q_{T}}\left( \frac{\widetilde{u}_{t}^{2}}{4}+\left( L\widetilde{%
u}\right) ^{2}\right) \varphi _{\lambda ,s}^{2}dxdt+ 
\]%
\begin{equation}
+C_{1}\lambda s\int\limits_{Q_{T}}\left( \nabla \widetilde{u}\right)
^{2}\varphi _{\lambda ,s}^{2}dxdt+C_{1}\lambda ^{2}s^{2}\int\limits_{Q_{T}}%
\widetilde{u}^{2}\varphi _{\lambda ,s}^{2}\leq  \label{6.22}
\end{equation}%
\[
\leq C_{1}\int\limits_{Q_{T}}\left( \left\vert \nabla \widetilde{m}%
\right\vert ^{2}+\widetilde{m}^{2}\right) \varphi _{\lambda
,s}^{2}dxdt+C_{1}e^{3\lambda \left( T+2\right) ^{s}}\delta ^{2},\text{ }%
\forall \lambda \geq 1,\text{ }\forall s\geq s_{1}. 
\]

We now proceed similarly with inequality (\ref{6.14}). When doing so, we use
(\ref{6.61}), (\ref{6.15}), (\ref{6.17}) and Theorem 5.2. Then we obtain the
following analog of (\ref{6.22})%
\[
\sqrt{s}\int\limits_{Q_{Y}}\left( \nabla \widetilde{m}\right) ^{2}\varphi
_{\lambda ,s}^{2}dxdt+\lambda s^{2}\int\limits_{Q_{T}}\widetilde{m}%
^{2}\varphi _{\lambda ,s}^{2}dxdt\leq 
\]%
\begin{equation}
\leq C_{1}\int\limits_{Q_{T}}\left[ \left( L\widetilde{u}\right)
^{2}+\left\vert \nabla \widetilde{u}\right\vert ^{2}+\widetilde{u}^{2}\right]
\varphi _{\lambda ,s}^{2}dxdt+  \label{6.23}
\end{equation}%
\[
+C_{1}e^{3^{s}\lambda }\left\Vert \widetilde{m}\left( x,0\right) \right\Vert
_{H^{1}\left( \Omega \right) }^{2}+C_{1}e^{3\lambda \left( T+2\right)
^{s}}\delta ^{2},\text{ }\forall s\geq s_{1}. 
\]

By (\ref{6.7}), (\ref{6.9}) and triangle inequality $\left\Vert \widetilde{m}%
\left( x,0\right) \right\Vert _{H^{1}\left( \Omega \right) }\leq C_{1}.$
Hence, using (\ref{6.23}), we obtain%
\[
\int\limits_{Q_{Y}}\left[ \left( \nabla \widetilde{m}\right) ^{2}+\widetilde{%
m}^{2}\right] \varphi _{\lambda ,s}^{2}dxdt\leq 
\]%
\begin{equation}
\leq \frac{C_{1}}{\sqrt{s}}\int\limits_{Q_{T}}\left[ \left( L\widetilde{u}%
\right) ^{2}+\left\vert \nabla \widetilde{u}\right\vert ^{2}+\widetilde{u}%
^{2}\right] \varphi _{\lambda ,s}^{2}dxdt+  \label{6.24}
\end{equation}%
\[
+C_{1}e^{3^{s}\lambda }+C_{1}e^{3\lambda \left( T+2\right) ^{s}}\delta ^{2},%
\text{ }\forall s\geq s_{1}. 
\]%
Substituting (\ref{6.24}) in (\ref{6.22}), we obtain%
\[
\int\limits_{Q_{T}}\left( \frac{\widetilde{u}_{t}^{2}}{4}+\left( L\widetilde{%
u}\right) ^{2}\right) \varphi _{\lambda ,s}^{2}dxdt+ 
\]%
\[
+C_{1}\lambda s\int\limits_{Q_{T}}\left( \nabla \widetilde{u}\right)
^{2}\varphi _{\lambda ,s}^{2}dxdt+C_{1}s^{2}\int\limits_{Q_{T}}\widetilde{u}%
^{2}\varphi _{\lambda ,s}^{2}\leq 
\]%
\[
\leq \frac{C_{1}}{\sqrt{s}}\int\limits_{Q_{T}}\left[ \left( L\widetilde{u}%
\right) ^{2}+\left\vert \nabla \widetilde{u}\right\vert ^{2}+\widetilde{u}%
^{2}\right] \varphi _{\lambda ,s}^{2}dxdt+ 
\]%
\[
+C_{1}e^{3^{s}\lambda }+C_{1}e^{3\lambda \left( T+2\right) ^{s}}\delta ^{2},%
\text{ }\forall \lambda \geq 1,\text{ }\forall s\geq s_{1}, 
\]%
The left hand side of this estimate dominates the term in its third line.
Hence,%
\[
\int\limits_{Q_{T}}\left( \frac{\widetilde{u}_{t}^{2}}{4}+\left( L\widetilde{%
u}\right) ^{2}\right) \varphi _{\lambda ,s}^{2}dxdt+\int\limits_{Q_{T}}\left[
\left( \nabla \widetilde{u}\right) ^{2}+\widetilde{u}^{2}\right] \varphi
_{\lambda ,s}^{2}dxdt\leq 
\]%
\begin{equation}
\leq C_{1}e^{3^{s}\lambda }+C_{1}e^{3\lambda \left( T+2\right) ^{s}}\delta
^{2},\text{ }\forall s\geq s_{1}.  \label{6.25}
\end{equation}%
Combining (\ref{6.25}) with (\ref{6.24}), we obtain%
\begin{equation}
\int\limits_{Q_{T}}\left( \frac{\widetilde{u}_{t}^{2}}{4}+\left( L\widetilde{%
u}\right) ^{2}+\left( \nabla \widetilde{u}\right) ^{2}+\widetilde{u}%
^{2}+\left( \nabla \widetilde{m}\right) ^{2}+\widetilde{m}^{2}\right)
\varphi _{\lambda ,s}^{2}dxdt\leq  \label{6.26}
\end{equation}%
\[
\leq C_{1}e^{3^{s}\lambda }+C_{1}e^{3\lambda \left( T+2\right) ^{s}}\delta
^{2},\text{ }\forall \lambda \geq 1,\text{ }\forall s\geq s_{1}. 
\]

By (\ref{5.11}) and (\ref{6.12})%
\[
\varphi _{\lambda ,s}\left( t\right) =e^{2\lambda \left( t+2\right)
^{s}}\geq e^{2\lambda \left( \gamma +2\right) ^{s}}\text{ in }Q_{\gamma T}. 
\]%
Hence, the first line of (\ref{6.26}) can be estimated from the below as: 
\[
\int\limits_{Q_{T}}\left( \frac{\widetilde{u}_{t}^{2}}{4}+\left( L\widetilde{%
u}\right) ^{2}+\left( \nabla \widetilde{u}\right) ^{2}+\widetilde{u}%
^{2}+\left( \nabla \widetilde{m}\right) ^{2}+\widetilde{m}^{2}\right)
\varphi _{\lambda ,s}^{2}dxdt\geq 
\]%
\[
\geq e^{2\lambda \left( \gamma +2\right) ^{s}}\int\limits_{Q_{T}}\left( 
\frac{\widetilde{u}_{t}^{2}}{4}+\left( L\widetilde{u}\right) ^{2}+\left(
\nabla \widetilde{u}\right) ^{2}+\widetilde{u}^{2}+\left( \nabla \widetilde{m%
}\right) ^{2}+\widetilde{m}^{2}\right) dxdt. 
\]%
Substituting this in (\ref{6.26}) and recalling (\ref{6.19}), we obtain 
\begin{equation}
\int\limits_{Q_{\gamma T}}\left( \frac{\widetilde{u}_{t}^{2}}{4}+\left( L%
\widetilde{u}\right) ^{2}+\left( \nabla \widetilde{u}\right) ^{2}+\widetilde{%
u}^{2}+\left( \nabla \widetilde{m}\right) ^{2}+\widetilde{m}^{2}\right)
dxdt\leq  \label{6.27}
\end{equation}%
\[
\leq C_{1}\exp \left[ -2\lambda \left( \gamma +2\right) ^{s}\left( 1-\frac{1%
}{2}\left( \frac{3}{\gamma +2}\right) \right) ^{s}\right] + 
\]%
\[
+C_{1}e^{3\lambda \left( T+2\right) ^{s}}\delta ^{2},\text{ }\forall \lambda
\geq 1,\forall s\geq s_{1}. 
\]%
Since by (\ref{6.122}) 
\[
1-\frac{1}{2}\left( \frac{3}{\gamma +2}\right) ^{s}>\frac{1}{2},\text{ }%
\forall s>0, 
\]%
then (\ref{6.27}) implies%
\begin{equation}
\int\limits_{Q_{\gamma T}}\left( \frac{\widetilde{u}_{t}^{2}}{4}+\left( L%
\widetilde{u}\right) ^{2}+\left( \nabla \widetilde{u}\right) ^{2}+\widetilde{%
u}^{2}+\left( \nabla \widetilde{m}\right) ^{2}+\widetilde{m}^{2}\right)
dxdt\leq  \label{6.28}
\end{equation}%
\[
\leq C_{1}e^{-\lambda \left( \gamma +2\right) ^{s}}+C_{1}e^{3\lambda \left(
T+2\right) ^{s}}\delta ^{2},\text{ }\forall \lambda \geq 1,\forall s\geq
s_{1}. 
\]%
Set $s=s_{1}$ and choose $\lambda =\lambda \left( \delta \right) $ such that 
\begin{equation}
e^{3\lambda \left( T+2\right) ^{s_{1}}}\delta ^{2}=\delta .  \label{6.29}
\end{equation}%
Hence, 
\begin{equation}
\lambda =\ln \left( \delta ^{-\alpha }\right) ,\text{ }\alpha =\frac{1}{%
3\left( T+2\right) ^{s_{1}}}.  \label{6.30}
\end{equation}%
Hence, to ensure that $\lambda \geq 1$ as in (\ref{6.19}), we should have $%
\delta \in \left( 0,\delta _{0}\right) ,$ where $\delta _{0}=\delta
_{0}\left( T,M_{1},M,\gamma \right) $ is so small that%
\begin{equation}
\delta _{0}\in \left( 0,e^{-1/\alpha }\right) .  \label{6.31}
\end{equation}%
Thus, (\ref{6.27})-(\ref{6.31}) lead to:%
\begin{equation}
\left. 
\begin{array}{c}
\int\limits_{Q_{\gamma T}}\left( \widetilde{u}_{t}^{2}/4+\left( L\widetilde{u%
}\right) ^{2}+\left( \nabla \widetilde{u}\right) ^{2}+\widetilde{u}%
^{2}+\left( \nabla \widetilde{m}\right) ^{2}+\widetilde{m}^{2}\right)
dxdt\leq \\ 
\leq C_{1}\delta ^{2\beta },\text{ }\forall \delta \in \left( 0,\delta
_{0}\right) ,%
\end{array}%
\right.  \label{6.32}
\end{equation}%
\begin{equation}
2\beta =\frac{1}{3}\left( \frac{\gamma +2}{T+2}\right) ^{s_{1}}\in \left( 0,%
\frac{1}{3}\right) .  \label{6.33}
\end{equation}%
The target estimate (\ref{6.123}) of this theorem follows immediately from
estimates (\ref{6.32}) and (\ref{6.33}). \ 

To prove uniqueness, set in (\ref{6.60}), (\ref{6.61}) $\delta =0.$ Then, we
obtain instead of (\ref{6.28}) 
\[
\int\limits_{Q_{\gamma T}}\left( \frac{\widetilde{u}_{t}^{2}}{4}+\left( L%
\widetilde{u}\right) ^{2}+\left( \nabla \widetilde{u}\right) ^{2}+\widetilde{%
u}^{2}+\left( \nabla \widetilde{m}\right) ^{2}+\widetilde{m}^{2}\right)
dxdt\leq 
\]%
\[
\leq C_{1}e^{-\lambda \left( \gamma +2\right) ^{s}},\text{ }\forall \lambda
\geq 1,\forall s\geq s_{1},\text{ }\forall \gamma \in \left( 0,T\right) . 
\]%
Setting here $s=s_{1},\lambda \rightarrow \infty ,$ we obtain $\widetilde{u}=%
\widetilde{m}=0$ in $Q_{\gamma T},$ $\forall \gamma \in \left( 0,T\right) .$ 
$\square $

\section{A Specification of Theorem 6.1 for the Retrospective Problem of
Section 4}

\label{sec:7}

In our specific Retrospective Problem of section 4, the principal parts of
elliptic operators of both equations (\ref{4.4}) and (\ref{4.5}) of the MFGS
are the same:%
\begin{equation}
L_{0}=\frac{\sigma _{1}^{2}(\mathbf{x})}{2}\partial _{x}^{2}+\frac{\sigma
_{2}^{2}(\mathbf{x})}{2}\partial _{y}^{2},  \label{7.1}
\end{equation}%
which fits well the generalized MFGS (\ref{6.3}), (\ref{6.4}). However,
unlike (\ref{6.4}), the analog of the term $F_{2}\left( \nabla m,m,u,\nabla
u,Lu,x,t\right) $ in (\ref{4.5}) contains the following terms with $u_{xx}$
and $u_{yy}:$%
\[
\frac{\varphi _{1}^{2}(\mathbf{x})}{a_{0}(\mathbf{x})}mu_{xx}+\frac{\varphi
_{2}^{2}(\mathbf{x})}{b_{0}(\mathbf{x})}mu_{yy}. 
\]%
Hence, we need to modify for this case both the formulation and the proof of
Theorem 6.1 via replacing the term $\left( L\widetilde{u}\right) ^{2}$ with
a stronger term $\widetilde{u}_{xx}^{2}+\widetilde{u}_{xy}^{2}+\widetilde{u}%
_{yy}^{2}.$ This can be done via a new Carleman estimate, which is a
specification of the Carleman estimate of Theorem 5.1. Below the domain $%
\Omega $ is the square in (\ref{2.2}). Consider the following modification
of the space $H_{0}^{2,1}\left( Q_{T}\right) $ in (\ref{5.10}): 
\[
\widehat{H}_{0}^{2,1}\left( Q_{T}\right) =\left\{ u\in H^{2,1}\left(
Q_{T}\right) :u\mid _{\Gamma _{0,T}}=0,\partial _{\nu }u\mid _{\Gamma
_{1,T}}=0\right\} . 
\]

\textbf{Theorem 7.1} (Carleman estimate: a specification of Theorem 5.1 for
Retrospective Problem of section 4).\emph{\ Assume that conditions (\ref%
{4.90})-(\ref{4.902}) hold. Let }$L_{0}$\emph{\ be the elliptic operator
defined in (\ref{7.1}). Let }$\varphi _{\lambda ,s}\left( t\right) $\emph{\
be the CWF (\ref{5.11}). There exists a number }$\widehat{C}=\widehat{C}%
\left( T,\sigma _{0},D\right) >0$\emph{\ and a sufficiently large number }$%
s_{2}=s_{2}\left( T,D\right) >1$, \emph{both numbers} \emph{depending only
on listed parameters, such that the following Carleman estimate holds}%
\[
\int\limits_{Q_{T}}\left( u_{t}+L_{0}u\right) ^{2}\varphi _{\lambda
,s}^{2}dxdt\geq \frac{1}{4}\int\limits_{Q_{T}}u_{t}^{2}\varphi _{\lambda
,s}^{2}dxdt+\widehat{C}\int\limits_{Q_{T}}\left(
u_{xx}^{2}+u_{xy}^{2}+u_{yy}^{2}\right) \varphi _{\lambda ,s}^{2}dxdt+ 
\]%
\begin{equation}
+\widehat{C}\lambda s\int\limits_{Q_{T}}\left( \nabla u\right) ^{2}\left(
t+2\right) ^{s-1}\varphi _{\lambda ,s}^{2}dxdt+\frac{1}{2}\lambda
^{2}s^{2}\int\limits_{Q_{T}}\left( t+2\right) ^{2s-2}u^{2}\varphi _{\lambda
,s}^{2}-  \label{7.2}
\end{equation}%
\[
-\lambda s\left( T+2\right) ^{s-1}e^{2\lambda \left( T+2\right)
^{s}}\int\limits_{\Omega }u^{2}\left( x,T\right) dx-\sigma
_{0}^{2}e^{2\lambda \left( T+2\right) ^{2}}\int\limits_{\Omega }\left(
\nabla u\left( x,T\right) \right) ^{2}dx, 
\]%
\[
\forall u\in \widehat{H}_{0}^{2,1}\left( Q_{T}\right) ,\text{ }\forall
\lambda \geq 1,\forall s\geq s_{2}. 
\]

\textbf{Proof.} Below in this section $\widehat{C}=\widehat{C}\left(
T,\sigma _{0},D\right) >0$ denotes different numbers depending only on
listed parameters. To prove this theorem, we only need to estimate from the
below the term with $\left( Lu\right) ^{2}$ in (\ref{5.12}) when $L$ is
replaced with $L_{0}.$ It is convenient to assume in this proof that $u\in
C^{3}\left( \overline{Q}_{T}\right) \cap \widehat{H}_{0}^{2,1}\left(
Q_{T}\right) .$ Next, density arguments lead to $u\in \widehat{H}%
_{0}^{2,1}\left( Q_{T}\right) .$

We have 
\[
\left( L_{0}u\right) ^{2}\varphi _{\lambda ,s}^{2}\left( t\right) =\left( 
\frac{\sigma _{1}^{2}(\mathbf{x})}{2}u_{xx}+\frac{\sigma _{2}^{2}(\mathbf{x})%
}{2}u_{yy}\right) ^{2}\varphi _{\lambda ,s}^{2}\left( t\right) = 
\]%
\[
=\left[ \left( \frac{\sigma _{1}^{2}(\mathbf{x})}{2}\right) ^{2}u_{xx}^{2}+%
\frac{\sigma _{1}^{2}(\mathbf{x})\sigma _{2}^{2}(\mathbf{x})}{2}%
u_{xx}u_{yy}+\left( \frac{\sigma _{2}^{2}(\mathbf{x})}{2}\right)
^{2}u_{yy}^{2}\right] \varphi _{\lambda ,s}^{2}\left( t\right) \geq 
\]%
\[
\geq \sigma _{0}^{2}\left( u_{xx}^{2}+u_{yy}^{2}\right) \varphi _{\lambda
,s}^{2}\left( t\right) + 
\]%
\begin{equation}
+\left( \frac{\sigma _{1}^{2}(\mathbf{x})\sigma _{2}^{2}(\mathbf{x})}{2}%
u_{xx}u_{y}\varphi _{\lambda ,s}^{2}\left( t\right) \right) _{y}-\frac{%
\sigma _{1}^{2}(\mathbf{x})\sigma _{2}^{2}(\mathbf{x})}{2}%
u_{xyx}u_{y}\varphi _{\lambda ,s}^{2}\left( t\right) =  \label{7.3}
\end{equation}%
\[
=\sigma _{0}^{2}\left( u_{xx}^{2}+u_{yy}^{2}\right) \varphi _{\lambda
,s}^{2}\left( t\right) +\left( \frac{\sigma _{1}^{2}(\mathbf{x})\sigma
_{2}^{2}(\mathbf{x})}{2}u_{xx}u_{y}\varphi _{\lambda ,s}^{2}\left( t\right)
\right) _{y}+ 
\]%
\[
+\left( -\frac{\sigma _{1}^{2}(\mathbf{x})\sigma _{2}^{2}(\mathbf{x})}{2}%
u_{xy}u_{y}\varphi _{\lambda ,s}^{2}\left( t\right) \right) _{x}+\frac{%
\sigma _{1}^{2}(\mathbf{x})\sigma _{2}^{2}(\mathbf{x})}{2}u_{xy}^{2}\varphi
_{\lambda ,s}^{2}\left( t\right) + 
\]%
\[
+\left( \frac{\sigma _{1}^{2}(\mathbf{x})\sigma _{2}^{2}(\mathbf{x})}{2}%
\right) _{x}u_{xy}u_{y}\varphi _{\lambda ,s}^{2}\left( t\right) . 
\]%
By Cauchy-Schwarz inequality, (\ref{4.90}) and (\ref{4.902})%
\[
\frac{\sigma _{1}^{2}(\mathbf{x})\sigma _{2}^{2}(\mathbf{x})}{2}%
u_{xy}^{2}\varphi _{\lambda ,s}^{2}\left( t\right) +\left( \frac{\sigma
_{1}^{2}(\mathbf{x})\sigma _{2}^{2}(\mathbf{x})}{2}\right)
_{x}u_{xy}u_{y}\varphi _{\lambda ,s}^{2}\left( t\right) \geq 
\]%
\[
\geq \sigma _{0}^{2}u_{xy}^{2}\varphi _{\lambda ,s}^{2}\left( t\right) -%
\widehat{C}\left( \nabla u\right) ^{2}\varphi _{\lambda ,s}^{2}\left(
t\right) . 
\]%
Combining this with (\ref{7.3}) and integrating over $Q_{T},$ we obtain%
\[
\int\limits_{Q_{T}}\left( L_{0}u\right) ^{2}\varphi _{\lambda ,s}^{2}\left(
t\right) dxdt\geq \sigma _{0}^{2}\int\limits_{Q_{T}}\left(
u_{xx}^{2}+u_{xy}^{2}+u_{yy}^{2}\right) \varphi _{\lambda ,s}^{2}\left(
t\right) dxdt- 
\]%
\begin{equation}
-\widehat{C}\int\limits_{Q_{T}}\left( \nabla u\right) ^{2}\varphi _{\lambda
,s}^{2}\left( t\right) .  \label{7.4}
\end{equation}%
Combining (\ref{7.4}) with (\ref{5.12}) and noticing that the term in the
second line of (\ref{7.4}) will be absorbed by the term 
\[
C\lambda s\int\limits_{Q_{T}}\left( \nabla u\right) ^{2}\left( t+2\right)
^{s-1}\varphi _{\lambda ,s}^{2}dxdt, 
\]%
we obtain (\ref{7.2}). $\square $

Because of (\ref{4.5}), (\ref{4.50}) and (\ref{4.51}), we replace the set $%
B\left( M_{1}\right) $ in (\ref{6.9})\emph{\ }with the set $Z\left(
M_{1}\right) ,$ which we define as%
\[
Z\left( M_{1}\right) =\left\{ 
\begin{array}{c}
\left( u,m\right) \in C^{2,1}\left( \overline{Q}_{T}\right) : \\ 
\max_{t\in \left[ 0,T\right] }\left\Vert u\left( \mathbf{x},t\right)
\right\Vert _{C^{2}\left( \overline{\Omega }\right) }\leq M_{1}, \\ 
\max_{t\in \left[ 0,T\right] }\left\Vert m\left( \mathbf{x},t\right)
\right\Vert _{C^{1}\left( \overline{\Omega }\right) }\leq M_{1}, \\ 
m\left( \mathbf{x},t\right) \geq 0\text{ in }Q_{T}%
\end{array}%
\right\} . 
\]

\textbf{Theorem 7.2} (a specification of Theorem 6.1 for Retrospective
Problem of section 4).\emph{\ Let (\ref{4.09})-(\ref{4.11}) and (\ref{6.12})
hold. Also, let the term with the function }$g$\emph{\ in (\ref{4.4}) has
the form (\ref{4.50}), where the function }$g$\emph{\ satisfies condition (%
\ref{2.4}). Let two pairs of functions }$\left( u_{1},m_{1}\right) ,\left(
u_{2},m_{2}\right) \in Z\left( M_{1}\right) $\emph{\ satisfy equations (\ref%
{4.4}), (\ref{4.5}) and boundary conditions (\ref{4.6}), (\ref{4.9}).
Suppose that (\ref{4.13}) and (\ref{4.14}) hold. Let }$T>1$\emph{\ and let (%
\ref{6.122}) be valid.} \emph{Then there exists a sufficiently small number }%
$\delta _{1}=\delta _{1}\left( T,M_{1},\gamma ,\sigma _{0},D,\varepsilon
\right) \in \left( 0,1\right) $\emph{\ and a number }$\rho =\rho \left(
T,M_{1},\gamma ,\sigma _{0},D,\varepsilon \right) \in \left( 0,1/6\right) ,$%
\emph{\ both these numbers depend only on listed parameters, such that the
following H\"{o}lder stability estimate for the Retrospective Problem of
section 4 is valid:}%
\[
\left\Vert \widetilde{u}_{t}\right\Vert _{L_{2}\left( Q_{\gamma T}\right)
}+\left\Vert \widetilde{u}_{xx}\right\Vert _{L_{2}\left( Q_{\gamma T}\right)
}+\left\Vert \widetilde{u}_{xy}\right\Vert _{L_{2}\left( Q_{\gamma T}\right)
}+\left\Vert \widetilde{u}_{yy}\right\Vert _{L_{2}\left( Q_{\gamma T}\right)
}+ 
\]%
\[
+\left\Vert \nabla \widetilde{u}\right\Vert _{L_{2}\left( Q_{\gamma
T}\right) }+\left\Vert \widetilde{u}\right\Vert _{L_{2}\left( Q_{\gamma
T}\right) }+ 
\]%
\[
+\left\Vert \nabla \widetilde{m}\right\Vert _{L_{2}\left( Q_{\gamma
T}\right) }+\left\Vert \widetilde{m}\right\Vert _{L_{2}\left( Q_{\gamma
T}\right) }\leq C_{2}\delta ^{\rho },\text{ }\forall \delta \in \left(
0,\delta _{1}\right) , 
\]%
\emph{where the number }$C_{2}=C_{2}\left( T,M_{1},\gamma ,\sigma
_{0},D,\varepsilon \right) >0$\emph{\ depends only on listed parameters. In
addition, problem (\ref{4.4})-(\ref{4.6}), (\ref{4.9}) has at most one
solution }$\left( u,m\right) \in B\left( M_{1}\right) .$

\textbf{Proof}. In this proof, $C_{2}=C_{2}\left( T,M_{1},\gamma ,\sigma
_{0},D,\varepsilon \right) >0$ denotes different numbers depending only on
listed parameters. We now follow the proof of Theorem 6.1, although we
replace (\ref{6.11}) with (\ref{4.52}). Using (\ref{7.2}), we obtain instead
of (\ref{6.22}) 
\[
\frac{1}{4}\int\limits_{Q_{T}}\widetilde{u}_{t}^{2}\varphi _{\lambda
,s}^{2}dxdt+\widehat{C}\int\limits_{Q_{T}}\left( \widetilde{u}_{xx}^{2}+%
\widetilde{u}_{xy}^{2}+\widetilde{u}_{yy}^{2}\right) \varphi _{\lambda
,s}^{2}dxdt+ 
\]%
\begin{equation}
+C_{2}\lambda s\int\limits_{Q_{T}}\left( \nabla \widetilde{u}\right)
^{2}\varphi _{\lambda ,s}^{2}dxdt+C_{2}\lambda ^{2}s^{2}\int\limits_{Q_{T}}%
\widetilde{u}^{2}\varphi _{\lambda ,s}^{2}\leq  \label{7.7}
\end{equation}%
\[
\leq C_{2}\int\limits_{Q_{T}}\left( \left\vert \nabla \widetilde{m}%
\right\vert ^{2}+\widetilde{m}^{2}\right) \varphi
_{1,s}^{2}dxdt+C_{1}e^{3\lambda \left( T+2\right) ^{s}}\delta ^{2},\text{ }%
\forall \lambda \geq 1,\forall s\geq s_{2}. 
\]%
Next, using Theorem 5.2, (\ref{4.5}), (\ref{4.91}) and (\ref{4.92}), we
obtain the following analog of (\ref{6.24}) 
\[
\int\limits_{Q_{Y}}\left[ \left( \nabla \widetilde{m}\right) ^{2}+\widetilde{%
m}^{2}\right] \varphi _{\lambda ,s}^{2}dxdt\leq 
\]%
\begin{equation}
\leq \frac{C_{2}}{\sqrt{s}}\int\limits_{Q_{T}}\left[ \widetilde{u}_{xx}^{2}+%
\widetilde{u}_{yy}^{2}+\left\vert \nabla \widetilde{u}\right\vert ^{2}+%
\widetilde{u}^{2}\right] \varphi _{\lambda ,s}^{2}dxdt+  \label{7.8}
\end{equation}%
\[
+C_{2}e^{3^{s}\lambda }+C_{2}e^{3\lambda \left( T+2\right) ^{s}}\delta ^{2},%
\text{ }\forall \lambda \geq 1,\text{ }\forall s\geq s_{2}. 
\]%
Given (\ref{7.7}) and (\ref{7.8}), the rest of the proof is the same as the
rest of the proof of Theorem 6.1 after (\ref{6.24}). $\square $

\begin{center}
\textbf{Acknowledgments}
\end{center}

The effort of the first author was supported by the US National Science
Foundation grant DMS 2436227. The effort of the second author was supported
by the Russian Science Foundation grant N24-21-00031.


\begin{thebibliography}{99}
\bibitem{A} Y.~Achdou, P.~Cardaliaguet, F.~Delarue, A.~Porretta, and
F.~Santambrogio. \newblock {\em Mean Field Games}, volume 2281 of \emph{%
Lecture Notes in Mathematics, C.I.M.E. Foundation Subseries}. \newblock %
Springer Nature, Cetraro, Italy, 2019,
https://doi.org/10.1007/978-3-030-59837-2.

\bibitem{Bauso} D.~Bauso, H.~Tembine, and T.~Basar. \newblock Opinion
dynamics in social networks through mean-field games. 
\newblock {\em SIAM J.
Control Optim}, 54:3225--3257, 2016, https://doi.org/10.1137/140985676.

\bibitem{Carm} R. Carmona, F. Delarue, Probabilistic Theory of Mean Field
Games with Applications, Vol. I, Springer, 2018,
https://doi.org/10.1007/978-3-319-58920-6.

\bibitem{Chow} Y.~T. Chow, S.~W. Fung, S.~Liu, L.~Nurbekyan, and S.~Osher. 
\newblock {A numerical algorithm for inverse problem from partial boundary
	measurement arising from mean field game problem}. 
\newblock {\em Inverse
Problems}, 39(1):014001, 12 2022, https://doi.org/10.1088/1361-6420/aca5b0.

\bibitem{Ding} M.~H. Ding, H. Liu and G.~H. Zheng, Determining a stationary
mean field game system from full/partial boundary measurement, SIAM J.
Mathematical Analysis, 57, 661--681, 2025,
https://doi.org/10.1137/23M1594327.

\bibitem{Huang1} M.~Huang, R.~P. Malham\'{e}, and P.~E. Caines. 
\newblock {Large population stochastic dynamic games: closed-loop McKean-Vlasov
	systems and the Nash certainty equivalence principle}. 
\newblock {\em
Commun. Inf. Syst.}, 6:221--251, 2006,
https://dx.doi.org/10.4310/CIS.2006.v6.n3.a5.

\bibitem{Huang2} M.~Huang, P.~E. Caines, and R.~P. Malham\'{e}. 
\newblock {Large-population cost-coupled LQG problems with nonuniform agents:
	individual-mass behavior and decentralized Nash equilibria}. 
\newblock {\em
IEEE Trans. Automat. Control}, 52:1560--1571, 2007,
http://dx.doi.org/10.1109/TAC.2007.904450.

\bibitem{KL} M.V. Klibanov and J. Li, \emph{Inverse Problems and Carleman
Estimates: Global Uniqueness, Global Convergence and Experimental Data}, De
Gruyter, 2021.

\bibitem{MFG1} M.~V. Klibanov and Y.~Averboukh. 
\newblock {Lipschitz stability estimate and uniqueness in the retrospective
	analysis for the mean field games system via two Carleman estimates}. %
\newblock {\em SIAM J. Mathematical Analysis}, 56:616--636, 2024,
https://doi.org/10.1137/23M1554801.

\bibitem{MFG2} M.~V. Klibanov, The mean field games system: Carleman
estimates, Lipschitz stability and uniquenes, J. Inverse and Ill-Posed
Problems, 31, 455-466, 2024, http://dx.doi.org/10.1515/jiip-2023-0023.

\bibitem{MFG6} M.~V. Klibanov, A coefficient inverse problem for the mean
field games system, \emph{Applied Mathematics and Optimization}, 88:54,
2023, http://dx.doi.org/10.1007/s00245-023-10042-0.

\bibitem{MFG3} M.~V. Klibanov, J.~ Li and H.~Liu, On the mean field games
system with lateral Cauchy data via Carleman estimates, J. Inverse and
Ill-Posed Problems, 32, 277-295, 2024,
http://dx.doi.org/10.1515/jiip-2023-0089.

\bibitem{MFG4} M.~V. Klibanov, J. ~Li and H. ~Liu, H\"{o}lder stability and
uniqueness for the mean field games system via Carleman estimates, Studies
in Applied Mathematics, 151, 1447--1470, 2023,
http://dx.doi.org/10.1111/sapm.12633.

\bibitem{MFG7} M.~V. Klibanov, J. ~Li \ and Z. ~Yang, Convexification
numerical method for the retrospective problem of mean field games, Applied
Mathematics and Optimization, 90:6, 2024,
https://doi.org/10.1007/s00245-024-10152-3.

\bibitem{MFGIPI} M.~V. Klibanov, J. ~Li \ and Z. ~Yang, Convexification
numerical method for a coefficient inverse problem for the system of
nonlinear parabolic equations governing mean field games, Inverse Problems
and Imaging, 19, 219-252, 2025, http://dx.doi.org/10.3934/ipi.2024031.

\bibitem{MFGbook} M.~V. Klibanov and J. ~Li, \emph{Carleman Estimates in
Mean Field Games}, De Gruyter, 2025.

\bibitem{Kol} V.~N. Kolokoltsov and O.~A. Malafeyev. 
\newblock {\em Many Agent Games in Socio-Economic Systems: Corruption,
	Inspection, Coalition Building, Network Growth, Security}. \newblock %
Springer Nature Switzerland AG, 2019,
https://doi.org/10.1007/978-3-030-12371-0.

\bibitem{Lad} O. ~A. Ladyzhenskaya, Boundary Value Problems of Mathematical
Physics, Springer, Berlin, 1985, https://doi.org/10.1007/978-1-4757-4317-3.

\bibitem{LL1} J.~M. Lasry and P.~L. Lions, \newblock Jeux \`{a} champ moyen.
i. le cas stationnaire. \newblock {\em C. R. Math. Acad. Sci. Paris},
343:619--625, 2006, https://doi.org/10.1016/j.crma.2006.09.019.

\bibitem{LL2} J.~M. Lasry and P.~L. Lions, \newblock Mean field games. %
\newblock {\em Japanese Journal of Mathematics}, 2:229--260, 2007,
http://dx.doi.org/10.1007/s11537-007-0657-8.

\bibitem{LRS} M.~M. Lavrentiev, V.~G. Romanov and S.~P. Shishatskii, %
\newblock Ill-Posed Problems of Mathematical Physics and Analysis,\newblock %
AMS, Providence, R.I., 1986,

\bibitem{Liao} Z. Liao and Q. L\"{u}, On inverse problems for mean field
games with common noise via Carleman estimate, Inverse Problems, 41 (2025)
045009, http://dx.doi.org/10.1088/1361-6420/adc3b9.

\bibitem{Liu1} H. Liu, C. Mou and S. Zhang, Inverse problems for mean field
games, Inverse Problems, 39, 085003, 2023,
https://doi.org/10.1088/1361-6420/acdd90.

\bibitem{Liu2} H. Liu and C.W.K. Lo, Determining state space anomalies in
mean field games, Nonlinearity, 38, 025010, 2025,
https://doi.org/10.1088/1361-6544/ada67d.

\bibitem{Ren1} K. Ren, N. Soedjak, K. Wang and H. Zhai, Reconstructing a
state-independent cost function in a mean-field game model, Inverse
Problems, 40, 105010, 2024, https://doi.org/10.1088/1361-6420/ad7497.

\bibitem{Ren2} K. Ren, N. Soedjak and K. Wang, Unique determination of cost
functions in a multipopulation mean field game model, J. Differential
Equations, 427, 2025, 843-867, https://doi.org/10.1016/j.jde.2025.02.037.

\bibitem{Trusov} N.~V. Trusov. \newblock Numerical study of the stock market
crises based on mean field games approach. 
\newblock {\em J. Inverse
Ill-Posed Probl.}, 29:849--865, 2021, https://doi.org/10.1515/jiip-2020-0016.

\bibitem{V} M. ~M. Vajnberg, Variational Method and Method of Monotone
Operators in the Theory of Nonlinear Equations, Israel Program for
Scientific Translations, Jerusalem, 1973.
\end{thebibliography}
\end{document}